\documentclass[reqno,a4paper]{amsart}
\usepackage{graphicx}
\usepackage{amsmath}
\usepackage{mathrsfs}
\usepackage{amsfonts}
\usepackage{amssymb}
\usepackage{enumerate}
\usepackage{latexsym}
\usepackage{graphpap}
\usepackage{cite}
\usepackage{tikz}
\usepackage{multirow}
\usepackage{makecell}

\newtheorem{theorem}{Theorem}[section]
\newtheorem{corollary}[theorem]{Corollary}
\newtheorem{lemma}[theorem]{Lemma}
\newtheorem{proposition}[theorem]{Proposition}

\theoremstyle{definition}

\newtheorem{example}[theorem]{Example}

\newtheorem{remark}[theorem]{Remark}

\newtheorem*{algo}{Algorithm}

\theoremstyle{remark}
\numberwithin{equation}{section}

\newcommand{\bfa}{\mathbf{a}}
\newcommand{\bfb}{\mathbf{b}}

\newcommand{\bfl}{\mathbf{l}}
\newcommand{\bfp}{\mathbf{p}}
\newcommand{\bfq}{\mathbf{q}}
\newcommand{\bfr}{\mathbf{r}}
\newcommand{\bfs}{\mathbf{s}}
\newcommand{\bft}{\mathbf{t}}
\newcommand{\bfu}{\mathbf{u}}
\newcommand{\bfv}{\mathbf{v}}
\newcommand{\bfw}{\mathbf{w}}

\newcommand{\precw}{\prec_\bfw}

\newcommand{\ocs}{\operatorname{\mathsf{occset}}}

\newcommand*{\hypo}{{\mathsf{hypo}}}

\newcommand{\ev}{\operatorname{\mathsf{ev}}}
\newcommand{\inver}{\operatorname{\mathsf{inv}}}

\newcommand{\op}{^*}
\newcommand{\var}{\mathsf{Var}}
\newcommand{\fb}{finitely based}
\newcommand{\nfb}{non-finitely based}
\newcommand{\inmon}{involution monoid}
\newcommand{\insem}{involution semigroup}

\newcommand{\con}{\operatorname{\mathsf{con}}}
\newcommand{\lin}{\operatorname{\mathsf{lin}}}
\newcommand{\mix}{\operatorname{\mathsf{mix}}}
\newcommand{\ml}{\operatorname{\mathsf{ml}}}

\newcommand{\occ}{\operatorname{\mathsf{occ}}}
\newcommand{\chaos}{\operatorname{\mathsf{chaos}}}

\begin{document}
\begin{sloppypar}

\title[Representations and identities of $(\hypo_n,^\sharp)$]{Representations and identities of hypoplactic monoids with involution}%

\thanks{This research was partially supported by the National Natural Science Foundation of China (Nos.
12271224, 12171213, 12161062), the Fundamental Research Funds for the Central University (No. lzujbky-2023-16) and the Natural Science Foundation of Gansu Province (No. 23JRRA1055).}

\author[B. B. Han]{Bin Bin Han$^{1}$}
\author[W. T. Zhang]{Wen Ting Zhang$^{1,\star}$}\thanks{$^\star$ Corresponding author} %
\author[Y. F. Luo]{Yan Feng Luo$^{1}$}
\author[J. X. Zhao]{Jin Xing Zhao$^{2}$}

\address{$^{1}$ School of Mathematics and Statistics, Lanzhou University, Lanzhou, Gansu 730000, PR China} %
\email{hanbb19@lzu.edu.cn, zhangwt@lzu.edu.cn, luoyf@lzu.edu.cn}

\address{$^{2}$ School of Mathematical Sciences, Inner Mongolia University, Hohhot, Inner Mongolia 010021, PR China}
\email{zhjxing@imu.edu.cn}



\begin{abstract}
Let $(\hypo_n,^\sharp)$ be the hypoplactic monoid of finite rank $n$ with Sch\"{u}tzenberger's involution $^{\sharp}$.
In this paper, we exhibit a faithful representation of $(\hypo_n,^\sharp)$ as an involution monoid of upper triangular matrices over any semiring
from a large class including the tropical semiring under the skew transposition. We then give a transparent combinatorial characterization of the word identities satisfied by $(\hypo_n,^\sharp)$. Further, we prove that $(\hypo_n,^\sharp)$ is non-finitely based if and only if $n=2, 3$ and give a polynomial time algorithm to check whether a given word identity holds in $(\hypo_n,^\sharp)$.
\end{abstract}

\maketitle

\section{Introduction}
The plactic monoid, a kind of tableaux algebra introduced by  Knuth \cite{K71} in the 1970s, was studied in detail by Lascoux and Sch\"{u}tzenberger \cite{LS81}. The plactic monoid arises from the combinatorics of Young tableaux by identifying words over a fixed ordered alphabet whenever they produce the same tableau via Schensted's insertion algorithm \cite{Sch61}. Plactic-like monoids are also tableaux algebras which arise from the combinatorics of tableaux as the plactic monoid, including the hypoplactic monoid \cite{KT97,Nov00}, the stalactic monoid \cite{HNT07,Pri13}, the taiga monoid \cite{Pri13}, the sylvester monoid \cite{HNT05}, the baxter monoid \cite{Gir12} and the stylic monoid\cite{AR22}. These monoids have attracted much attention due to their interesting connection with combinatorics \cite{Lot02} and applications in symmetric functions \cite{Mac08}, representation theory \cite{Ful97}, Kostka-Foulkes polynomials \cite{LS78,LS81} and Schubert polynomials \cite{LS85, LS90}.

Identities and varieties of algebras have long been studied, and several important questions arise in this area.
An \textit{identity basis} for an algebra $A$ is a set of identities satisfied by $A$ that axiomatize all identities satisfied by $A$. An algebra $A$ is said to be \textit{finitely based} (FB) if it has some finite identity basis. Otherwise, it is said to be \textit{non-finitely based} (NFB). The \textit{finite basis problem} for algebras, that is the problem of classifying algebras according to the finite basis property, is one of the most prominent research problems in universal algebra. It is known that finite groups \cite{OP64}, finite associative rings \cite{Kru73,Lvov73}, and finite lattices \cite{Mck70} are
finitely based, but not all finite algebras are finitely based. In the 1960s, Perkins gave the first example of a {\nfb} finite semigroup ~\cite{Perkins66, Perkins69}, since then the finite basis problem for semigroups has attracted much attention. Refer to the survey \cite{Vol01} for more information on the finite basis problem for finite semigroups.
For a given semigroup $S$, the \textit{identity checking problem} \textsc{Check-Id}($S$) is the decision problem whose
instance is an arbitrary identity $\bfu \approx \bfv$, and the answer to such an instance is `YES'
if $S$ satisfies $\bfu \approx \bfv$, and `NO' if it does not.
For a finite semigroup $S$, the identity checking problem \textsc{Check-Id}($S$) is always decidable,  but this is not necessarily true for an infinite semigroup \cite{Mu68}.
Studying the computational complexity of identity checking in semigroups and other `classical' algebras such as groups and rings was proposed
by Sapir \cite[Problem~2.4]{KS95}, and many results in this area have been obtained so far\cite{AVG09,JM06,DJK18,Chen2020, KV20}.

Recall that a unary operation~$^*$ on a semigroup~$S$ is an \textit{involution} if~$S$ satisfies the identities
\begin{equation}
(x^*)^* \approx  x \quad \text{and} \quad (xy)^* \approx y^*x^*. \label{id: inv}
\end{equation}
An \textit{involution semigroup} is a pair $(S, \op)$ where~$S$ is a semigroup with {involution}~$^*$, and
$S$ is called the \textit{semigroup reduct} of $(S,\op)$.  Common examples of involution semigroups include groups $(G,^{-1})$ with inversion~$^{-1}$, multiplicative matrix semigroups $(M_n,^T)$ over any field with transposition~$^T$ and multiplicative matrix semigroups $(M_n,^D)$ over any field with skew transposition~$^D$.
Over the years, the identities and varieties of {\insem}s have received less attention than those of semigroups. However,
since the turn of the millennium, interest in {\insem}s has significantly increased.
For example, many counterintuitive results were established and examples have been found
of involution semigroups $(S,^*)$ and their semigroup reducts $S$ which are not simultaneously finitely based \cite{GZL-TB2,GZL-A01,Lee16a,Lee16b,Lee17a,Lee18,Lee19}. Refer to \cite{Lee20,Lee19+,GZL-variety, ADPV14, ADV12a, ADV12b} for more information on the identities and varieties of involution semigroups.

Recently the representations and identities of plactic monoids and  plactic-like monoids have received a lot of attention.
Johnson et al. have given a faithful representation of the plactic  monoid of each finite rank as a monoid of upper triangular matrices over the tropical semiring \cite[Theorem~2.9]{JK19}. Since the monoid of all upper triangular $n\times n$  tropical matrices satisfies non-trivial identities \cite{IZ14}, every plactic monoid of finite rank satisfies non-trivial identities. However, the plactic monoid of infinite rank does not satisfy any non-trivial
identity \cite[Theorem~3.2]{CKK17}, and so the plactic monoid of infinite rank is finitely based. The plactic monoid of rank $2$
is non-finitely based by \cite[Remark 4.6]{JK19}, and \cite[Corollary 5.4]{CHLS16} or
alternatively \cite[Theorem 1]{Shne89} and \cite[Theorem 4.1]{DJK18}, while the
plactic monoid of rank $3$ is non-finitely based by \cite[Corollary 4.5]{JK19} and
\cite[Theorem 3.15]{HZL}. The plactic monoid of rank $3$ with involution is non-finitely based by \cite[Remark 3.17]{HZL}, a similar argument with  \cite[Remark 3.17]{HZL} can show that the plactic monoid of rank $2$ with involution is non-finitely based. The finite basis problem for the plactic monoids of rank greater than or equal to $4$ and their involution cases are still open. Cain et al. have given  faithful representations of the hypoplactic, stalactic, taiga, sylvester and baxter monoids of each finite rank as monoids of upper triangular matrices over any semiring from a large class including the tropical semiring and fields of characteristic $0$ \cite{CJKM21}.  Cain and Malheiro have shown that the hypoplactic,
stalactic, taiga, sylvester and baxter monoids satisfy non-trivial identities \cite{CM18}.
In \cite{CMR21a} [resp. \cite{CMR21b, CJKM21,HZ}], it is shown that all hypoplactic [resp. stalactic, taiga, sylvester and baxter] monoids of rank greater than or equal to $2$ generate the same variety and are finitely based.
Aird and Ribeiro have given a faithful representation of the stylic monoid of each finite rank as a monoid of upper unitriangular matrices over tropical semiring, and then solved the finite basis problems for the stylic monoid and its involution case \cite{TD22}.
Volkov has solved the finite basis problem for the stylic monoid by different means \cite{volkov2022}.
The identity checking problems for the sylvster, baxter and stylic monoids have been shown to be decidable in polynomial time \cite{CMR21b,TD22}.
By the characterization of the identities satisfied by the hypoplactic, stalactic and taiga monoids \cite{CMR21b, CJKM21,HZ}, it is easy to show that the identity checking problems for them are in the complexity class $\mathsf{P}$.

The hypoplactic monoid $\hypo_n$ of finite rank $n$, first studied in depth by Novelli \cite{Nov00}, can be obtained by factoring the free monoid $\mathcal{A}_n^{\star}=\mathcal{A}_n^{+}\cup \{\varepsilon\}$ over the finite ordered alphabet $\mathcal{A}_n = \{1 < 2 < \cdots < n\}$ by a congruence that can be defined by the Krob-Thibon insertion algorithm that computes a combinatorial object from a word. Its elements can be uniquely identified with quasi-ribbon tableaux, and it is presented by the Knuth relations and the hypoplactic relations $cadb\equiv acbd$ for $a \leq b < c \leq d$ and $bdac\equiv dbca$ for $a < b \leq c < d$ with $a, b, c, d \in \mathcal{A}_n$. The hypoplactic monoid $\hypo_n$ forms an involution monoid, denoted by $(\hypo_n, ^{\sharp})$, under the the anti-automorphism $^{\sharp}$ induced by the unique order-reversing permutation on $\mathcal{A}_n$, which is called \textit{Sch\"{u}tzenberger's involution}.

Cain et al. have given faithful representations of $\hypo_n$ for each finite $n$ \cite[Theorem~3.4]{CJKM21} and proved that $\hypo_n$ for all $n\geq 2$ satisfy exactly the same identities and are finitely based \cite[Theorems~3.7 and 4.8]{CMR21a}.
However, it is easy to check that the representation  of $\hypo_n$ given in \cite[Theorem~3.4]{CJKM21} can not be extended to the involution case.

In this paper, we investigate the representations and identities of $(\hypo_n,^\sharp)$ for all finite $n$.
In Section \ref{sec:repre}, we exhibit a faithful representation of $(\hypo_n,^\sharp)$ for each finite $n$ as an involution monoid of upper triangular matrices over any semiring from a large class including the tropical semiring under the skew transposition. Clearly, our representation of $(\hypo_n,^\sharp)$ is also a representation of $\hypo_n$. Further, we prove that $(\hypo_n,^\sharp)$ for all $n\geq 4$ satisfy exactly the same identities.
In Section \ref{sec:characterization}, we give a characterization of the word identities satisfied by $(\hypo_n,^\sharp)$.
In Section \ref{sec:FBP4}, we prove that $(\hypo_2,^\sharp)$ and $(\hypo_3,^\sharp)$ are non-finitely based and $(\hypo_n,^\sharp)$ with  $n\geq 4$ is finitely based; we also show that each variety generated by $(\hypo_n,^\sharp)$ with $n\geq 3$ contains continuum many subvarieties.
Finally, we give a polynomial time algorithm to check whether a given word identity holds in $(\hypo_n,^\sharp)$ in Section \ref{sec:CHECK-ID}.
It is easy to see that the presence of $^\sharp$ changes the representations and identities of $\hypo_n$ quite radically, we summarize some properties of $\hypo_n$ and $(\hypo_n,^\sharp)$ in Table~\ref{table:compare}.

\begin{table*}[t]
\centering
\caption{Some properties of $\hypo_n$ and $(\hypo_n,^\sharp)$}
\begin{tabular}{|c|c|c|c|c|c|c|c|c|c|c|c|}
\hline
Properties & {\makecell[c]{$\hypo_n$\\ with $n\geq 2$}} & $(\hypo_2,^\sharp)$ & $(\hypo_3,^\sharp)$ &{\makecell[l]{$(\hypo_n,^\sharp)$\\ with $n\geq 4$}} \\
\hline
{\makecell[c]{Matrix \\representations}} & {\makecell[l]{Theorem\,3.4\\ in \cite{CJKM21}}} & Theorem\,3.1 & Theorem\,3.2 & Theorem\,3.6 \\
\hline
{\makecell[c]{Characterization \\of its identities}}  & {\makecell[l]{Theorem\,4.1\\ in \cite{CMR21a}}} & Theorem\,4.3 & Theorem\,4.5 & Theorem\,4.7 \\
\hline
{\makecell[c]{Finite basis \\property}}	& FB & NFB	& NFB & FB \\
\hline
{\makecell[c]{Number of \\subvarieties of \\its variety}}& $2^{\aleph_0}$ & unknown & $2^{\aleph_0}$ & $2^{\aleph_0}$ \\
\hline
Axiomatic rank 	& 4	& infinite	& infinite	& 4	\\
\hline
{\makecell[c]{Identity \\checking problem}}   &$\mathsf{P}$ & $\mathsf{P}$ & $\mathsf{P}$ & $\mathsf{P}$ \\
\hline
\end{tabular}
\label{table:compare}
\end{table*}

\section{Preliminaries} \label{sec: prelim}
Most of the notation and definitions of this article are given in this section.
Refer to the monograph of Burris and Sankappanavar~\cite{BS81} for any undefined notation and terminology of universal algebra in general.

\subsection{Words}

Let~$\mathcal{X}$ be a nonempty alphabet and let $\mathcal{X}^* =\{x^* \,|\, x \in \mathcal{X}\}$ be a disjoint copy of~$\mathcal{X}$.
Elements of $\mathcal{X} \cup \mathcal{X}^*$ are called \textit{variables}, elements of the free involution monoid $F_{\mathsf{inv}}^{\varepsilon}(\mathcal{X})=(\mathcal{X} \cup \mathcal{X}^*)^+ \cup \{\varepsilon\}$ are called \textit{words}, and elements of the free monoid $\mathcal{X}^{\star}=\mathcal{X}^+ \cup \{\varepsilon\}$ are called \textit{plain words}. A word $\bfu$ is a \textit{factor} of a word $\bfw$ if $\bfw = \bfa\bfu\bfb$ for some $\bfa, \bfb \in F_{\mathsf{inv}}^{\varepsilon}(\mathcal{X})$.

Let $\bfu\in (\mathcal{X} \cup \mathcal{X}^*)^+$ be a word and $x \in \mathcal{X} \cup \mathcal{X}^*$ be a variable. The \textit{content} of $\bfu$, denoted by $\con (\bfu)$, is the set of variables that occur in~$\bfu$, the \textit{length} of $\bfu$ is the number of variables occurring in $\bfu$ and is denoted by $|\bfu|$, and $\occ(x, \bfu)$ is the number of occurrences of $x$ in $\bfu$. Let $\overline{\bfu}$ be the plain word obtained from $\bfu$ by removing all occurrences of the symbol $^*$.
For $x_1, x_2, \ldots, x_n \in \mathcal{X}\cup\mathcal{X}^*$ such that $\overline{x_1}, \overline{x_2}, \ldots, \overline{x_n}\in \mathcal{X}$ are distinct variables, let $\bfu[x_1, x_2, \ldots, x_n]$ denote the word obtained from~$\bfu$ by retaining only the variables $x_1, x_1^*, x_2, x_2^*, \ldots, x_n, x_n^*$.

\begin{example}
Let $\bfu = x^* zxy^*x$ for some $x,y,z\in \mathcal{X}$. Then
\begin{itemize}
  \item  $\con (\bfu)=\{x, z, x^*, y^*\}$, $|\bfu|=5$, $\overline{\bfu}=xzxyx$;
  \item $\occ (x, \bfu) = 2$, $\occ (x^*,\bfu) = \occ (z,\bfu) =\occ (y^*,\bfu)=1$;
  \item $\bfu[x]=x^*x^2, \bfu[x, y]=x^*xy^*x$.
\end{itemize}
\end{example}

We use ${_i}x$ to refer to the $i$-th from the left occurrence of $x$ in $\bfu$ and ${_\infty}x$ to refer to the last occurrence of $x$ in $\bfu$.  The set $\ocs(\bfu) = \{ {_i}x \mid x \in \mathcal{X} \cup \mathcal{X}^*, 1 \le i \le \occ(x, \bfu) \}$ of all occurrences of all variables in $\bfu$ is called the {\em occurrence set of $\bfu$}. The word $\bfu$ induces an order $\prec_\bfu$ on the set $\ocs(\bfu)$ defined by ${_i}x \prec_\bfu {_j}y$ if and only if the $i$-th occurrence of $x$ precedes the $j$-th occurrence of $y$ in $\bfu$. We write
\[
\{{_{m_1}}x_1, {_{m_2}}x_2, \dots,  {_{m_r}}x_r\} \prec_\bfu \{{_{n_1}}y_1, {_{n_2}}y_2, \dots, {_{n_t}}y_t\}
\]
to mean that ${_{m_i}}x_i \prec_\bfu {_{n_j}}y_j$ for all $i$ and $j$. In particular, if $r=1$ [resp. $t=1$], we write ${_{m_1}}x_1  \prec_\bfu \{{_{n_1}}y_1, {_{n_2}}y_2, \dots, {_{n_t}}y_t\}$ [resp. $\{{_{m_1}}x_1, {_{m_2}}x_2, \dots, {_{m_r}}x_r\} \prec_\bfu {_{n_1}}y_1$]
rather than $\{{_{m_1}}x_1\}  \prec_\bfu \{{_{n_1}}y_1, {_{n_2}}y_2, \dots, {_{n_t}}y_t\}$ [resp. $\{{_{m_1}}x_1, {_{m_2}}x_2, \dots, {_{m_r}}x_r\} \prec_\bfu \{{_{n_1}}y_1\}$]. For convenience, if $m_1=m_2=\dots=m_r=\infty$ and  $n_1=n_2=\dots=n_t=1$, we also write $\{x_1, x_2, \dots, x_r\} \prec_\bfu \{y_1, y_2, \dots, y_t\}$.

\subsection{Terms and identities}

The set $\mathsf{T}(\mathcal{X})$ of \textit{terms} over~$\mathcal{X}$ is the smallest set containing $\mathcal{X}$ that is closed under concatenation and $^*$.
The proper inclusion $F_{\mathsf{inv}}^{\varepsilon}(\mathcal{X}) \subset \mathsf{T}(\mathcal{X})$ holds and the identities \eqref{id: inv} can be used to convert any nonempty term $\bft \in \mathsf{T}(\mathcal{X})$ into some unique word $\lfloor \bft\rfloor \in (\mathcal{X} \cup \mathcal{X}^*)^+$.
For instance, $\lfloor  x^2(x(yx^*)^*)^* zy^*\rfloor = x^2y(x^*)^2zy^*$.

An \textit{identity} is an expression $\bfs \approx \bft$ formed by nonempty terms $\bfs, \bft \in \mathsf{T}(\mathcal{X})$, a \textit{word identity} [resp. \textit{plain identity}] is an identity $\bfu \approx \bfv$ formed by words $\bfu, \bfv \in (\mathcal{X} \cup \mathcal{X}^*)^+$ [resp. $\mathcal{X}^+$]. The identity $\overline{\bfu}\approx \overline{\bfv}$ is called the \textit{plain projection} of $\bfu\approx \bfv$.
We write $\bfu = \bfv$ if $\bfu$ and $\bfv$ are identical.  An identity $\bfu \approx \bfv$ is \textit{non-trivial} if
$\bfu \neq \bfv$.
An  {\insem} $(S,\op)$ \textit{satisfies} an identity $\bfs \approx \bft$, if for any substitution $\varphi: \mathcal{X} \to S$, the elements $\bfs\varphi$ and $\bft\varphi$ of $S$ coincide; in this case, $\bfs \approx \bft$ is also said to be an \textit{identity of}  $(S,\op)$.

Clearly any  {\inmon} that satisfies a  word identity $\bfs \approx \bft$ also satisfies the identity $\bfs [x_1, x_2, \ldots, x_n] \approx \bft[x_1, x_2, \ldots, x_n]$ for any distinct variables $\overline{x_1}, \overline{x_2}, \ldots, \overline{x_n} \in \mathcal{X}$, since assigning the unit element~$1$ to a variable~$x$ in a  word identity is effectively the same as removing all occurrences of $x$ and~$x^*$.

For any semigroup $S$ [resp. {\insem} $(S,\op)$],
a set $\Sigma$ of identities of $S$ [resp. $(S,\op)$] is an \textit{identity basis} for $S$ [resp. $(S,\op)$] if every identity of $S$ [resp. $(S,\op)$] is deducible from $\Sigma$.
A semigroup [resp. {\insem}] is \textit{\fb} if it has some finite identity basis; otherwise, it is \textit{\nfb}.

The variety generated by a semigroup $S$ [resp. {\insem} $(S,\op)$] is denoted by $\var S$ [resp. $\var (S,\op)$]. For any set $\Sigma$ of identities, denote by $\var \Sigma$ the variety determined by $\Sigma$.

The \textit{axiomatic rank} of a semigroup $S$ [resp. {\insem} $(S,\op)$] is the least natural number $r$ such that $S$ [resp. $(S,\op)$] admits an identity basis $\Sigma$, where the number of occurrences of distinct variables in each identity [resp. the plain projection of each identity] in $\Sigma$ does not exceed $r$. If no such natural number exists, we say that $S$ [resp. $(S,\op)$] has \textit{infinite axiomatic rank}. Note that if $S$ [resp. $(S,\op)$] is finitely based, then it has finite axiomatic rank.

\subsection{The hypoplactic monoid and its involution}%
Let $\mathcal{A} = \{ 1 < 2 < 3 < \cdots \}$ denote the set of positive integers, viewed as an infinite ordered alphabet.
The tableau and insertion algorithm related to the hypoplactic monoid are given in the following.

A \textit{quasi-ribbon tableau} is a finite grid of cells, aligned so that the leftmost cell in each row is below the rightmost cell of the previous row, filled with variables from $\mathcal{A}$, such that the entries in each row are weakly increasing from left to right, and the entries in each column are strictly increasing from top to bottom.
The associated insertion algorithm is as follows:

\begin{algo}[Krob--Thibon algorithm]\label{algo:hypo}
Input a quasi-ribbon tableau $T$ and a variable $a \in \mathcal{A}$.
If there is no entry in $T$ that is less than or equal to $a$, output the tableau obtained by creating a new cell, labelled with $a$, and gluing $T$ by its top-leftmost entry to the bottom of this new cell; otherwise let $x$ be the right-most and bottom-most entry of $T$ that is less than or equal to $a$. Separate $T$ in two parts, such that one part is from the top left down to and including $x$. Put a new entry $a$ to the right of $x$ and glue the remaining part of $T$ (below and to the right of $x$) onto the bottom of the new entry $a$.
Output the resulting tableau.
\end{algo}

Let $w_1, \cdots, w_k\in \mathcal{A}$ and $\bfw=w_1\cdots w_k \in \mathcal{A}^{\star}$. Then the quasi-ribbon tableau ${\rm P}_{\hypo_{\infty}}(\bfw)$ of $\bfw$ is obtained as follows:  reading
$\bfw$ from left-to-right, one starts with an empty tableau and inserts each variable in $\bfw$ into a quasi-ribbon tableau according to the above Algorithm.
For example, ${\rm P}_{\hypo_{\infty}}(36131512665)$ is given as follows:
\begin{equation*}
\begin{tikzpicture}
[ampersand replacement=\&,row sep=-\pgflinewidth,column sep=-\pgflinewidth]
\matrix [nodes=draw]
{
\node {1}; \& \node{1}; \& \node{1};\& \node{2}; \& \& \& \& \&\\
 \& \&  \& \node{3};  \& \node{3}; \&  \node{5}; \& \node{5};\& \&\\
 \&  \& \&  \& \& \& \node{6};\& \node{6}; \& \node{6};\\
};
\end{tikzpicture}
\end{equation*}
Note that the same variable cannot appear in two different rows of a quasi-ribbon tableau.

Define the relation $\equiv_{\hypo_{\infty}}$ by
\[
\bfu \equiv_{\hypo_{\infty}} \bfv \quad  \text{if and only if} \quad {\rm P}_{\hypo_{\infty}}(\bfu) = {\rm P}_{\hypo_{\infty}}(\bfv)
\]
for any $\bfu,\bfv \in \mathcal{A}^{\star}$. The relation $\equiv_{\hypo_{\infty}}$ is a congruence on $\mathcal{A}^{\star}$. The hypoplactic monoid $\hypo_{\infty}$ is the factor monoid $ \mathcal{A}^{\star}/_{\equiv_{\hypo_{\infty}}}$. The rank-$n$ analogue $\hypo_{n}$ is the factor monoid $ \mathcal{A}^{\star}_n/_{\equiv_{\hypo_{\infty}}}$, where the relation $\equiv_{\hypo_{\infty}}$ is naturally restricted to $\mathcal{A}^{\star}_n\times\mathcal{A}^{\star}_n$ and $\mathcal{A}_n = \{1 < 2 < \cdots < n\}$ is the set of the first $n$ natural numbers viewed as a finite ordered alphabet.
It follows from the definition of $\equiv_{\hypo_{\infty}}$ that each element $[\bfu]_{\equiv_{\hypo_{\infty}}}$
of the factor monoid $\hypo_{\infty}$ can be identified with the
combinatorial object ${\rm P}_{\hypo_{\infty}}(\bfu)$.
Clearly $\hypo_1$ is a free monogenic monoid $\langle 1 \rangle=\{\varepsilon, 1, 1^2, 1^3,\ldots\}$ and so it is commutative.
Note that
\begin{align}\label{M-order}
\hypo_1 \subset \hypo_2 \subset \cdots \subset \hypo_i \subset\hypo_{i+1}\subset \cdots \subset \hypo_{\infty}.
\end{align}

For any $\bfu \in \mathcal{A}^{\star}$, the \textit{evaluation} of $\bfu$, denoted by $\ev(\bfu)$, is the infinite
tuple of non-negative integers, indexed by $\mathcal{A}$, whose $a$-th element is $\occ(a, \bfu)$, thus this tuple describes the number of each variable in $\mathcal{A}$ that appears in $\bfu$. It is immediate from the definition of the hypoplactic monoid that if $\bfu \equiv_{\hypo_\infty} \bfv$, then $\ev(\bfu) = \ev(\bfv)$, and hence it makes sense to define the evaluation of each element of the hypoplactic monoid to be the evaluation
of any word representing it.  Note that $\ev(\bfu) = \ev(\bfv)$ implies $\con(\bfu) = \con(\bfv)$.

Let $\bfu \in \mathcal{A}^{\star}$ and $\con \left( \bfu \right) = \{a_1 < \dots < a_k\}$ for some $k \in \mathbb{N}$. We say $\bfu$ has an \textit{$a_{i+1}$-$a_i$ inversion} if for $1 \leq i \leq k-1$ there is at least one occurrence of $a_{i+1}$ before the last occurrence of $a_i$ when reading $\bfu$ from left-to-right. Note that we only consider inversions of consecutive elements of the content of $\bfu$. Denote by $\inver(\bfu)$ the set of all inversions occurring in $\bfu$. For example, $\inver(36131512665)=\{3\textrm{-}2, 6\textrm{-}5\}$.
	
\begin{proposition}[{\cite[Subsection~4.2]{Nov00}}]\label{pro:inversion}
For any $\bfu, \bfv \in \mathcal{A}^{\star}$, $\bfu \equiv_{\hypo_\infty} \bfv$ if and only if $\ev(\bfu)=\ev(\bfv)$ and $\inver(\bfu)=\inver(\bfv)$.
\end{proposition}

The hypoplactic monoid  can also be defined by the presentation $\left\langle \mathcal{A} \mid \mathcal{R}_{\hypo_\infty} \right\rangle$, where
\begin{align*}
\mathcal{R}_{\hypo_\infty} =& \left\{ (acb,cab): a \leq b < c \right\} \cup \left\{ (bac,bca): a < b \leq c \right\}\\
 \cup & \left\{ (cadb,acbd): a \leq b < c \leq d \right\} \cup \left\{ (bdac,dbca): a < b \leq c < d \right\}.
\end{align*}

For each $n \in \mathbb{N}$, a presentation for the hypoplactic monoid of rank $n$ can be obtained by restricting generators and relations of the above presentation to the generators in $\mathcal{A}_n$. Note that these relations preserve the evaluation of words.

Consider the anti-automorphism $^\sharp$ of $\mathcal{A}_n^{\star}$ given by reversing the linear order on $\mathcal{A}_n$. By \cite[Theorem~5.4]{Nov00}, the relation $\equiv_{\hypo_{\infty}}$ is compatible with the unary operation $^{\sharp}$, that is, for any $\bfw, \bfw' \in \mathcal{A}_n^{\star}$, $\bfw \equiv_{\hypo_{\infty}} \bfw'$ if and only if $\bfw^{\sharp} \equiv_{\hypo_{\infty}} {(\bfw')}^{\sharp}$. Hence $(\hypo_n, ^{\sharp})$ is an involution monoid, and the involution $^{\sharp}$ is called \textit{Sch\"{u}tzenberger's involution}.

\begin{proposition}
Sch\"{u}tzenberger's involution is the only involution on the hypoplactic monoid $\hypo_n$ for $n\geq 1$.
\end{proposition}
\begin{proof}
Suppose that $^*$ is an involution operation on $\hypo_n$. Note that the relations in $\mathcal{R}_{\hypo_\infty}$ preserve the evaluation of words. Thus the involution of a generator in $\mathcal{A}_n$ is still a generator in $\mathcal{A}_n$.
Since $\hypo_1$ has only one generator $1$, we have $1^*=1$. Thus the involution on $\hypo_1$ is trivial.  For $\hypo_n$ with $n\geq 2$, let $a< b\leq n$. Then $(aba)^*=a^*b^*a^*\equiv_{\hypo_{\infty}} a^*a^*b^*=(baa)^*$ by $aba \equiv_{\hypo_{\infty}} baa \in \mathcal{R}_{\hypo_\infty}$. This implies $a^*b^*a^*\equiv_{\hypo_{\infty}} a^*a^*b^*\in \mathcal{R}_{\hypo_\infty}$, and so $b^*< a^*$.  Hence for any $a<b$, there must be $b^*< a^*$ under the involution $^*$. Therefore $^*$ induces the unique order-reversing permutation on $\mathcal{A}_n$.
\end{proof}

\subsection{Matrix representations over semirings}
Recall that $\mathbb{S}=(S, +, \cdot)$ is a \textit{commutative semiring} with additive identity $\mathbf{0}$ and multiplicative identity $\mathbf{1}$ if $S$ is a set equipped with two binary operations $+$ and $\cdot$ such that $(S, +)$ and $(S, \cdot)$ are  commutative monoids satisfying
\[
a(b+c)=a\cdot b+a\cdot c \quad\text{and}\quad \mathbf{0}\cdot a=\mathbf{0}
\]
for all $a,b,c \in S$. We say that $\mathbb{S}$ is \textit{idempotent} if $a+a=a$ for all $a\in S$.
An element $a\in S$ has \textit{infinite multiplicative order} if for any non-negative integers $i,j$, $a^{i}=a^{j}$ if and only if $i=j$.
In this paper, we always assume that $\mathbb{S}$ is a commutative and idempotent semiring with $\mathbf{0}, \mathbf{1}$ containing an element of infinite multiplicative order.
A common example of such a semiring is the tropical semiring
$\mathbb{T}= (\mathbb{R} \cup \{-\infty\}, \oplus, \otimes)$, which is the set $\mathbb{R}$ of real numbers together with minus infinity $-\infty$, with the addition and multiplication defined as follows
\[
a \oplus b = \max \{a, b\} \quad \mbox{and} \quad a \otimes b = a + b.
\]
In other words, the tropical sum of two numbers is their maximum and the tropical
product of two numbers is their sum, and
$-\infty$ is the additive identity and $0$ is the multiplicative identity. Note that, except for $-\infty$ and $0$, all other elements in $\mathbb{T}$ have infinite multiplicative order.

It is easy to see that the set of all $n \times n$  matrices with entries in $\mathbb{S}$ forms a semigroup under the matrix multiplication induced from the operations in $\mathbb{S}$. We denote this semigroup by $M_{n}(\mathbb{S})$. Let $UT_n(\mathbb{S})$ be subsemigroup of $M_n(\mathbb{S})$ of all upper triangular $n \times n$  matrices.
For any matrix $A\in M_n(\mathbb{S})$, denote by $A^{D}$ the matrix obtained by reflecting $A$ with respect to the secondary diagonal (from the top right to the bottom left corner), that is,  $(A^{D})_{ij}=A_{(n+1-j)(n+1-i)}$. It is easy to verify that this unary operation $^D$ is an involution on $M_{n}(\mathbb{S})$ which is called the \textit{skew transposition}.
A (linear) representation of a semigroup $S$ [resp. involution semigroup $(S,^*)$] is a
homomorphism $\rho : S \rightarrow M_n(\mathbb{S})$ [resp. $\rho : (S,^*) \rightarrow (M_n(\mathbb{S}),^{D})$ ]. The homomorphism $\rho$ is said to be \textit{faithful} if it is injective. Note that any involution semigroup representation $\rho : (S,^*) \rightarrow (M_n(\mathbb{S}),^{D})$ induces a semigroup representation $\rho : S \rightarrow M_n(\mathbb{S})$, but a semigroup representation $\rho : S \rightarrow M_n(\mathbb{S})$ does not necessarily induce an involution semigroup representation $\rho : (S,^*) \rightarrow (M_n(\mathbb{S}),^{D})$. The tropical semiring is of interest as a natural carrier for representations of semigroups. For example,  the bicyclic monoid $\mathcal{B}:=\langle a,b\mid ba=1\rangle$, which is ubiquitous in infinite semigroup theory, admits no faithful finite dimensional representations over any field; however it has a number of natural representations over the tropical semiring \cite{DJK18,Izhakian10}.

\section{matrix representations of $(\hypo_n, ^\sharp)$}\label{sec:repre}%
In this section, we exhibit a faithful matrix representation of $(\hypo_n, ^\sharp)$ for each finite $n$ as an involution monoid of upper triangular matrices over $\mathbb{S}$ under the skew transposition, and we prove that all involution semigroups $(\hypo_n, ^{\sharp})$ with $n\geq 4$ generate the same variety.

For convenience, denote by
 \begin{gather*}
 \mathrm{P}=\begin{pmatrix}
s &  \mathbf{0} \\
\mathbf{0} &  \mathbf{1}
\end{pmatrix},
\mathrm{Q}=\begin{pmatrix}
\mathbf{1} & \mathbf{0} \\
\mathbf{0} & s
\end{pmatrix},
 \mathrm{J}=\begin{pmatrix}
\mathbf{1} & \mathbf{0} & \mathbf{0} \\
\mathbf{0} & \mathbf{1}& \mathbf{1}\\
\mathbf{0} & \mathbf{0} & \mathbf{1}
\end{pmatrix},
\mathrm{K}=\begin{pmatrix}
\mathbf{1} &\mathbf{1} &  \mathbf{0}\\
\mathbf{0} & \mathbf{1}& \mathbf{0} \\
\mathbf{0} & \mathbf{0} & \mathbf{1}
\end{pmatrix}
\end{gather*}
where $s\in \mathbb{S}$ is an element of infinite multiplicative order. Denote by $\mathsf{diag}\{\Lambda_1,\Lambda_2,\dots, \Lambda_n\}$ the block diagonal matrix
\[
\begin{pmatrix}
\Lambda_1 & & & \\
 & \Lambda_2 & &\\
 & & \ddots &\\
  &&& \Lambda_n
\end{pmatrix}
\]
where $\Lambda_1,\Lambda_2,\dots, \Lambda_n$ are square matrices. And let $\mathrm{E}_n$ be the $n\times n$ matrix with $\mathbf{1}$s on the main diagonal and $\mathbf{0}$s elsewhere.

First we give a matrix representation of $(\hypo_1,^\sharp)$.
Define a map $\psi_1: \mathcal{A}_1\cup \{\varepsilon\}\rightarrow UT_2(\mathbb{S})$ given by $\varepsilon\mapsto \mathrm{E}_2$ and
$1\mapsto \mathrm{P}\mathrm{Q}$.
Clearly, the map $\psi_1$ induces a faithful representation of $(\hypo_1,^\sharp)$.

Next we consider a matrix representation of $(\hypo_2,^\sharp)$.
Define a map $\psi_2: \mathcal{A}_2\cup \{\varepsilon\}\rightarrow UT_5(\mathbb{S})$ given by $\varepsilon\mapsto \mathrm{E}_5,$
\[
1\mapsto \mathsf{diag}\{s,\mathrm{J}, \mathbf{1}\},\;\; 2\mapsto\mathsf{diag}\{\mathbf{1},\mathrm{K}, s\}.
\]
Clearly, $\psi_2$ can be extended to a homomorphism from $\mathcal{A}_2^{\star}$ to $UT_5(\mathbb{S})$. Note that
$\psi_2(1^{\sharp})= \psi_2(2)=\mathsf{diag}\{\mathbf{1},\mathrm{K}, a\}=(\psi_2(1))^{D}$ and $\psi_2(2^{\sharp})= \psi_2(1)=\mathsf{diag}\{a,\mathrm{J}, \mathbf{1}\}=(\psi_2(2))^{D}$.
Thus $\psi_2$  can be extended to a homomorphism $\psi_2: (\mathcal{A}_2^{\star}, ^{\sharp})\rightarrow (UT_5(\mathbb{S}),^{D})$.
In fact, $\psi_2$ induces a faithful representation of $(\hypo_2,^\sharp)$.

\begin{theorem}\label{thm:hypo2-repre}
The map $\psi_2:(\hypo_2,^\sharp)\rightarrow(UT_5(\mathbb{S}),^{D})$ is a faithful representation of $(\hypo_2,^\sharp)$.
\end{theorem}

\begin{proof}
Note that $\psi_2$ is a homomorphism from $(\mathcal{A}_2^{\star}, ^{\sharp})$ to $(UT_5(\mathbb{S}),^{D})$. Then to show that $\psi_2$ induces a homomorphism from $(\hypo_2,^\sharp)$ to $(UT_5(\mathbb{S}),^{D})$, we only need to show that for any $\bfu,\bfv\in \mathcal{A}_2^{\star}$, if $\bfu\equiv_{\hypo_\infty}\bfv$, then
$\psi_2(\bfu) = \psi_2(\bfv)$.
By the definition of $\psi_2$, it is easy to verify that for any $\bfw\in \mathcal{A}_2^{\star}$,
\begin{align*}\label{id:2}
\psi_2(\bfw)=
\left\{
  \begin{array}{ll}
\vspace{0.1cm}
    \mathrm{E}_5, & \hbox{if $\bfw=\varepsilon$,} \\
\vspace{0.1cm}
    \mathsf{diag}\{s^{\occ(1, \bfw)},\mathrm{J}, \mathbf{1}\}, & \hbox{if $\con(\bfw)=\{1\}$,} \\
\vspace{0.1cm}
    \mathsf{diag}\{\mathbf{1},\mathrm{K}, s^{\occ(2, \bfw)}\}, & \hbox{if $\con(\bfw)=\{2\}$,} \\
\vspace{0.1cm}
   \mathsf{diag}\{s^{\occ(1, \bfw)}, \mathrm{KJ}, s^{\occ(2, \bfw)}\}, & \hbox{if $\con(\bfw)=\{1,2\}$ and $2\textrm{-}1 \in \inver(\bfw)$,} \\
\vspace{0.1cm}
    \mathsf{diag}\{s^{\occ(1, \bfw)}, \mathrm{JK},s^{\occ(2, \bfw)}\}, & \hbox{if $\con(\bfw)=\{1,2\}$ and $2\textrm{-}1 \not\in \inver(\bfw)$.}
  \end{array}
\right.
\end{align*}
Since $\ev(\bfu)=\ev(\bfv), \inver(\bfu)=\inver(\bfv)$ by Proposition \ref{pro:inversion}, it is routine to show that $\psi_2(\bfu) = \psi_2(\bfv)$.

Suppose that $\bfu\not\equiv_{\hypo_\infty}\bfv$. Then $\ev(\bfu)\neq\ev(\bfv)$ or $\inver(\bfu)\neq\inver(\bfv)$ by Proposition \ref{pro:inversion}. By the definition of $\psi_2$ it is easy to see that $\psi_2(\bfu)\neq\psi_2(\bfv)$, a contradiction.
Hence $\psi_2$ is injective. Therefore $\psi_2:(\hypo_2, ^\sharp)\rightarrow(UT_5(\mathbb{S}),^{D})$ is a faithful representation of $(\hypo_2,^\sharp)$.
\end{proof}

Next we consider a matrix representation of $(\hypo_3,^\sharp)$.
Define a map $\psi_3: \mathcal{A}_3\cup \{\varepsilon\}\rightarrow UT_{13}(\mathbb{S})$ given by
$\varepsilon\mapsto \mathrm{E}_{13}$,
\begin{gather*}
1\mapsto
\mathsf{diag}\{\mathrm{P},\mathrm{J},\mathrm{J},\mathrm{E}_3,\mathrm{E}_2\},\;\;
2\mapsto
\mathsf{diag}\{\mathrm{Q},\mathrm{K},\mathrm{KJ},\mathrm{J},\mathrm{P}\},\;\;
3\mapsto
\mathsf{diag}\{\mathrm{E}_2,\mathrm{E}_3, \mathrm{K},\mathrm{K},\mathrm{Q}\}.
\end{gather*}
Clearly, $\psi_3$ can be extended to a homomorphism from $\mathcal{A}_3^{\star}$ to $UT_{13}(\mathbb{S})$. Note that $\psi_3(1^{\sharp})=\psi_3(3)=\mathsf{diag}\{\mathrm{E}_2,\mathrm{E}_3, \mathrm{K},\mathrm{K},\mathrm{Q}\}=(\psi_3(1))^{D}, \psi_3(2^{\sharp})=\psi_3(2)=\mathsf{diag}\{\mathrm{Q},\mathrm{K}$, $\mathrm{KJ},\mathrm{J},\mathrm{P}\}=(\psi_3(2))^{D}$ and
$\psi_3(3^{\sharp})=\psi_3(1)=\mathsf{diag}\{\mathrm{P},\mathrm{J},\mathrm{J},\mathrm{E}_3,\mathrm{E}_2\}=(\psi_3(1))^{D}$.
Thus $\psi_3$ can be extended to a homomorphism $\psi_3: (\mathcal{A}_3^{\star}, ^{\sharp})\rightarrow (UT_{13}(\mathbb{S}),^{D})$.
In fact, $\psi_3$ induces a faithful representation of  $(\hypo_3,^\sharp)$.

\begin{theorem}\label{thm:hypo3-repre}
The map $\psi_3:(\hypo_3,^\sharp)\rightarrow(UT_{13}(\mathbb{S}),^{D})$ is a faithful representation of $(\hypo_3,^\sharp)$.
\end{theorem}

\begin{proof}
Note that $\psi_3$ is a homomorphism from $(\mathcal{A}_3^{\star}, ^{\sharp})$ to $(UT_{13}(\mathbb{S}),^{D})$. Then to show that the map $\psi_3$ induces a homomorphism from $(\hypo_3,^\sharp)$ to $(UT_{13}(\mathbb{S}),^{D})$, we only need to show that for any $\bfu,\bfv\in \mathcal{A}_3^{\star}$, if $\bfu\equiv_{\hypo_\infty}\bfv$, then $\psi_3(\bfu) = \psi_3(\bfv)$.
By the definition of $\psi_3$, it is easy to verify that for any $\bfw\in \mathcal{A}_3^{+}$,
\[
\psi_3(\bfw)=\mathsf{diag}\{\Lambda_1,\Lambda_2,\Lambda_3,\Lambda_4,\Lambda_5\}
\]
where $\Lambda_1=\mathrm{P}^{\occ(1, \bfw)}\mathrm{Q}^{\occ(2, \bfw)}$, $\Lambda_5=\mathrm{P}^{\occ(2, \bfw)}\mathrm{Q}^{\occ(3, \bfw)}$ and
\begin{align*}\label{id:3}
&\Lambda_2= \left\{
  \begin{array}{ll}
\mathrm{E_3}, & \hbox{if $\con(\bfw)=\{3\}$,} \\ [0.1cm]
\mathrm{J}, & \hbox{if $\con(\bfw)=\{1\}$ or $\{1,3\}$,} \\ [0.1cm]
\mathrm{K}, & \hbox{if $\con(\bfw)=\{2\}$ or $\{2,3\}$,} \\ [0.1cm]
\mathrm{KJ}, & \hbox{if $\{1,2\}\subseteq\con(\bfw)$ and $2\textrm{-}1 \in \inver(\bfw)$,} \\ [0.1cm]
\mathrm{JK}, & \hbox{if $\{1,2\}\subseteq\con(\bfw)$ and $2\textrm{-}1 \not\in \inver(\bfw)$,}
\end{array}
\right.\\
&\Lambda_3= \left\{
\begin{array}{ll}
\mathrm{J}, & \hbox{if $\con(\bfw)=\{1\}$,} \\ [0.1cm]
\mathrm{K}, & \hbox{if $\con(\bfw)=\{3\}$,} \\ [0.1cm]
\mathrm{JK}, & \hbox{if $\con(\bfw)=\{1,3\}$ and $3\textrm{-}1 \not\in \inver(\bfw)$,} \\[0.1cm]
\mathrm{KJ}, & \hbox{if $\{2\} \subseteq \con(\bfw)$, or $\con(\bfw)=\{1,3\}$ and $3\textrm{-}1 \in \inver(\bfw)$,}
\end{array}
\right.\\
&\Lambda_4= \left\{
  \begin{array}{ll}
\mathrm{E_3}, & \hbox{if $\con(\bfw)=\{1\}$,} \\ [0.1cm]
\mathrm{J}, & \hbox{if $\con(\bfw)=\{2\}$ or $\{1,2\}$,} \\ [0.1cm]
\mathrm{K}, & \hbox{if $\con(\bfw)=\{3\}$ or $\{1,3\}$,} \\ [0.1cm]
\mathrm{KJ},  & \hbox{if $\{2,3\}\subseteq\con(\bfw)$ and $3\textrm{-}2 \in \inver(\bfw)$,} \\ [0.1cm]
\mathrm{JK}, & \hbox{if $\{2,3\}\subseteq\con(\bfw)$ and $3\textrm{-}2 \not\in \inver(\bfw)$.}
\end{array}
\right.
\end{align*}
Since $\ev(\bfu)=\ev(\bfv), \inver(\bfu)=\inver(\bfv)$ by Proposition \ref{pro:inversion}, it is routine to show that $\psi_3(\bfu) = \psi_3(\bfv)$.

Suppose that $\bfu\not\equiv_{\hypo_\infty}\bfv$. Then $\ev(\bfu)\neq\ev(\bfv)$ or $\inver(\bfu)\neq\inver(\bfv)$ by Proposition \ref{pro:inversion}. By the definition of $\psi_3$ it is easy to see that $\psi_3(\bfu)\neq\psi_3(\bfv)$, a contradiction.
Hence $\psi_3$ is injective. Therefore $\psi_3:(\hypo_3, ^\sharp)\rightarrow(UT_{13}(\mathbb{S}),^{D})$ is a faithful representation of $(\hypo_3,^\sharp)$.
\end{proof}

Now we consider a matrix representation of $(\hypo_n,^\sharp)$ for $n\geq 4$.
For any $([\bfu]_{\hypo_3}, [\bfv]_{\hypo_3})\in \hypo_3\times \hypo_3$, we can define an involution operation $^{\sharp}$ on $\hypo_3\times \hypo_3$ by  $([\bfu]_{\hypo_3}, [\bfv]_{\hypo_3})^{\sharp}= ([\bfv]_{\hypo_3}^{\sharp}, [\bfu]_{\hypo_3}^{\sharp})$. For any $i<j \in \mathcal{A}_n$ with $n\geq 4$, it follows from the definition of $^{\sharp}$ that $j^{\sharp}<i^{\sharp}$ and there is at most one $k\in\mathcal{A}_n$ satisfying $k=k^{\sharp}$ and $j-i= i^{\sharp}-j^{\sharp}$. So there are five cases about the order of $i,j, i^{\sharp}, j^{\sharp}$ in $\mathcal{A}_n$: $i^{\sharp}=j$,  $i<j\leq j^{\sharp}<i^{\sharp}$, $j^{\sharp}<i^{\sharp}\leq i<j$, $i<j^{\sharp}<j<i^{\sharp}$, $j^{\sharp}<i<i^{\sharp}<j$.
For any $i<j \in \mathcal{A}_n$ with $n\geq 4$, we can define a map $\varphi_{ij}$ from $\mathcal{A}_n^{\star}$ to $\hypo_3\times \hypo_3$ which can be determined by the following three cases according to the order of $i, i^{\sharp}, j, j^{\sharp}$ in $\mathcal{A}_n$.

{\bf Case 1.} $i^{\sharp} = j$. Define a map $\lambda: \mathcal{A}_n \rightarrow \hypo_3$ given by
\begin{align*}
k &\mapsto \begin{cases}
[1]_{\hypo_3} & \text{if}\ k = i,\\
[3]_{\hypo_3} & \text{if}\ k = j,\\
[31]_{\hypo_3} & \text{if}\ i < k < j,\\
\left[\varepsilon\right]_{\hypo_3} & \text{otherwise}.
\end{cases}
\end{align*}
Clearly, this map can be extended to a homomorphism $\lambda: \mathcal{A}_n^{\star} \rightarrow \hypo_3$.
Define a map $\lambda_{ij}: \mathcal{A}_n \rightarrow \hypo_3\times \hypo_3$ given by
\[
k \mapsto (\lambda(k), \lambda(k)).
\]
This map can be extended to a homomorphism $\lambda_{ij}: \mathcal{A}_n^{\star} \rightarrow \hypo_3\times \hypo_3$. Further, $\lambda_{ij}$
is also a homomorphism from $(\mathcal{A}_n^{\star},^{\sharp})$ to $(\hypo_3\times \hypo_3,^{\sharp})$. This is because for any $k \in \mathcal{A}_n$, $\lambda(k^{\sharp})=(\lambda(k))^{\sharp}$ which follows from
{\setlength{\arraycolsep}{0.5pt}
\[
\left\{
\begin{array}{llll}
\lambda(k^{\sharp})=[3]_{\hypo_3}&=(\lambda(k))^{\sharp} & \quad \text{if}\ k = i,\\
\lambda(k^{\sharp})=[31]_{\hypo_3}&=(\lambda(k))^{\sharp} & \quad\text{if}\ i < k < j,\\
\lambda(k^{\sharp})=[1]_{\hypo_3}&=(\lambda(k))^{\sharp}  & \quad\text{if}\ k = j,\\
\lambda(k^{\sharp})=[\varepsilon]_{\hypo_3}&=(\lambda(k))^{\sharp}   & \quad\text{otherwise}.
\end{array}\right.
\]}
Therefore for any $\bfw=k_1k_2\cdots k_n$,
\begin{align*}
\lambda_{ij}(\bfw^{\sharp})&= (\lambda(\bfw^{\sharp}), \lambda(\bfw^{\sharp}))\\
&=(\lambda(k_n^{\sharp})\cdots\lambda(k_1^{\sharp}), \lambda(k_n^{\sharp})\cdots\lambda(k_1^{\sharp}))\\
&=((\lambda(k_n))^{\sharp}\cdots(\lambda(k_1))^{\sharp}, (\lambda(k_n))^{\sharp}\cdots(\lambda(k_1))^{\sharp})\\
&=((\lambda(\bfw))^{\sharp}, (\lambda(\bfw))^{\sharp})\\
&=(\lambda_{ij}(\bfw))^{\sharp}.
\end{align*}

{\bf Case 2.} $i<j\leq j^{\sharp}<i^{\sharp}$ or $j^{\sharp}<i^{\sharp}\leq i<j$. For convenience, let $i_1=i, i_2=j, i_3=j^{\sharp}$ and $i_4=i^{\sharp}$ when $i<j\leq j^{\sharp}<i^{\sharp}$ and $i_1=j^{\sharp}, i_2=i^{\sharp},i_3=i$ and $i_4=j$ when $j^{\sharp}<i^{\sharp}\leq i<j$. Define maps $\theta_1: \mathcal{A}_n \rightarrow \hypo_3$ and $\theta_2: \mathcal{A}_n \rightarrow \hypo_3$ by
\begin{align*}
k \mapsto \begin{cases}
[1]_{\hypo_3} & \text{if}\ k = i_1,\\
[2]_{\hypo_3} & \text{if}\ k = i_2,\\
[21]_{\hypo_3} & \text{if}\ i_1 < k < i_2,\\
\left[\varepsilon\right]_{\hypo_3} & \text{otherwise,}
\end{cases}\;\;\; \mbox{and}\;\;\;
k \mapsto \begin{cases}
[2]_{\hypo_3} & \text{if}\ k = i_3,\\
[3]_{\hypo_3} & \text{if}\ k = i_4,\\
[32]_{\hypo_3} & \text{if}\ i_3 < k < i_4,\\
\left[\varepsilon\right]_{\hypo_3} & \text{otherwise,}
\end{cases}
\end{align*}
respectively.
Clearly, $\theta_1, \theta_2$ can be extended to homomorphisms from $\mathcal{A}_n^{\star}$ to $\hypo_3$ respectively.
Define a map $\theta_{ij}: \mathcal{A}_n \rightarrow \hypo_3\times \hypo_3$ by
\[
k \mapsto (\theta_1(k), \theta_2(k)).
\]
This map can be extended to a homomorphism $\theta_{ij}: \mathcal{A}_n^{\star} \rightarrow \hypo_3\times \hypo_3$. Further, $\theta_{ij}$ is also a homomorphism from $(\mathcal{A}_n^{\star},^{\sharp})$ to $(\hypo_3\times \hypo_3,^{\sharp})$. This is because for any $k \in \mathcal{A}_n$, $\theta_1(k^{\sharp})=(\theta_2(k))^{\sharp}, (\theta_1(k))^{\sharp}=\theta_2(k^{\sharp})$ which follows from
{\setlength{\arraycolsep}{0.5pt}
\begin{align*}
\left\{
\begin{array}{llllllll}
\theta_1(k^{\sharp})&=[\varepsilon]_{\hypo_3}&=(\theta_2(k))^{\sharp},\quad
(\theta_1(k))^{\sharp}&=[3]_{\hypo_3}&= \theta_2(k^{\sharp})&\quad \text{if}~k = i_1,\\
\theta_1(k^{\sharp})&=[\varepsilon]_{\hypo_3}&=(\theta_2(k))^{\sharp},\quad
(\theta_1(k))^{\sharp}&=[32]_{\hypo_3}&= \theta_2(k^{\sharp})&\quad\text{if}\ i_1 < k < i_2,\\
\theta_1(k^{\sharp})&=[\varepsilon]_{\hypo_3}&=(\theta_2(k))^{\sharp},\quad
(\theta_1(k))^{\sharp}&=[2]_{\hypo_3} &=\theta_2(k^{\sharp})&\quad \text{if}\ k= i_2, i_2\neq i_3,\\
\theta_1(k^{\sharp})&=[2]_{\hypo_3}&=(\theta_2(k))^{\sharp},\quad
(\theta_1(k))^{\sharp}&=[2]_{\hypo_3}&= \theta_2(k^{\sharp})&\quad \text{if}\ k = i_2=i_3,\\
\theta_1(k^{\sharp})&=[2]_{\hypo_3}&=(\theta_2(k))^{\sharp},\quad
(\theta_1(k))^{\sharp}&=[\varepsilon]_{\hypo_3}&= \theta_2(k^{\sharp})&\quad \text{if}\ k= i_3, i_2\neq i_3,\\
\theta_1(k^{\sharp})&=[21]_{\hypo_3}&=(\theta_2(k))^{\sharp},\quad
(\theta_1(k))^{\sharp}&=[\varepsilon]_{\hypo_3}&= \theta_2(k^{\sharp})&\quad \text{if}\  i_3<k<i_4,\\
\theta_1(k^{\sharp})&=[1]_{\hypo_3}&=(\theta_2(k))^{\sharp},\quad
(\theta_1(k))^{\sharp}&=[\varepsilon]_{\hypo_3}&= \theta_2(k^{\sharp})&\quad \text{if}\ k =i_4,\\
\theta_1(k^{\sharp})&=[\varepsilon]_{\hypo_3}&=(\theta_2(k))^{\sharp},\quad
(\theta_1(k))^{\sharp}&=[\varepsilon]_{\hypo_3}&= \theta_2(k^{\sharp})&\quad\text{otherwise}.
\end{array}\right.
\end{align*}}
Therefore for any $\bfw=k_1k_2\cdots k_n$,
\begin{align*}
\theta_{ij}(\bfw^{\sharp})&= (\theta_1(\bfw^{\sharp}), \theta_2(\bfw^{\sharp}))\\
&=(\theta_1(k_n^{\sharp})\cdots\theta_1(k_1^{\sharp}), \theta_2(k_n^{\sharp})\cdots\theta_2(k_1^{\sharp}))\\
&=((\theta_2(k_n))^{\sharp}\cdots(\theta_2(k_1))^{\sharp}, (\theta_1(k_n))^{\sharp}\cdots(\theta_1(k_1))^{\sharp})\\
&=((\theta_2(\bfw))^{\sharp}, (\theta_1(\bfw))^{\sharp})\\
&=(\theta_{ij}(\bfw))^{\sharp}.
\end{align*}

{\bf Case 3.} $i<j^{\sharp}<j<i^{\sharp}$ or $j^{\sharp}<i<i^{\sharp}<j$. For convenience, let $i_1=i, i_2=j^{\sharp}, i_3=j$ and $i_4=i^{\sharp}$ when $i<j^{\sharp}<j<i^{\sharp}$ and $i_1=j^{\sharp}, i_2=i,i_3=i^{\sharp}$ and $i_4=j$ when $j^{\sharp}<i<i^{\sharp}<j$. Define  maps $\kappa_1: \mathcal{A}_n \rightarrow \hypo_3$ and $\kappa_2: \mathcal{A}_n \rightarrow \hypo_3$ by
\begin{align*}
k \mapsto \begin{cases}
[1]_{\hypo_3} & \text{if}\ k = i_1,\\
[21]_{\hypo_3} & \text{if}\ i_1 < k < i_3,\\
[2]_{\hypo_3} & \text{if}\  k= i_3,\\
\left[\varepsilon\right]_{\hypo_3} & \text{otherwise},
\end{cases} \;\;\;
k \mapsto \begin{cases}
[2]_{\hypo_3} & \text{if}\  k = i_2,\\
[3]_{\hypo_3} & \text{if}\ k = i_4,\\
[32]_{\hypo_3} & \text{if}\ i_2 < k < i_4,\\
\left[\varepsilon\right]_{\hypo_3} & \text{otherwise}
\end{cases}
\end{align*}
respectively.
Clearly, $\kappa_1, \kappa_2$ can be extended to homomorphisms from $\mathcal{A}_n^{\star}$ to $\hypo_3$ respectively.
Define a map $\kappa_{ij}: \mathcal{A}_n \rightarrow \hypo_3\times \hypo_3$,
\[
k \mapsto (\kappa_1(k), \kappa_2(k)).
\]
Clearly, this map can be extended to a homomorphism from $\mathcal{A}_n^{\star}$ to $\hypo_3\times \hypo_3$. Further, $\kappa_{ij}$ is a homomorphism from $(\mathcal{A}_n^{\star},^{\sharp})$ to $(\hypo_3\times \hypo_3,^{\sharp})$.
This is because for any $k \in \mathcal{A}_n$, $\kappa_1(k^{\sharp})=(\kappa_2(k))^{\sharp}, (\kappa_1(k))^{\sharp}=\kappa_2(k^{\sharp})$ which follows from
{\setlength{\arraycolsep}{0.5pt}
\begin{align*}
\left\{
\begin{array}{lllllllll}
\kappa_1(k^{\sharp})&=[\varepsilon]_{\hypo_3}&=(\kappa_2(k))^{\sharp},\quad
(\kappa_1(k))^{\sharp}&=[3]_{\hypo_3}&=\kappa_2(k^{\sharp})&\quad \text{if}\ k = i_1,\\
\kappa_1(k^{\sharp})&=[\varepsilon]_{\hypo_3}&=(\kappa_2(k))^{\sharp},\quad
(\kappa_1(k))^{\sharp}&=[32]_{\hypo_3}&=\kappa_2(k^{\sharp})& \quad \text{if}\ i_1 < k < i_2,\\
\kappa_1(k^{\sharp})&=[2]_{\hypo_3}&=(\kappa_2(k))^{\sharp},\quad
(\kappa_1(k))^{\sharp}&=[32]_{\hypo_3}&=\kappa_2(k^{\sharp})& \quad\text{if}\ k = i_2,\\
\kappa_1(k^{\sharp})&=[21]_{\hypo_3}&=(\kappa_2(k))^{\sharp},\quad
(\kappa_1(k))^{\sharp}&=[32]_{\hypo_3}&= \kappa_2(k^{\sharp})& \quad\text{if}\ i_2< k< i_3,\\
\kappa_1(k^{\sharp})&=[21]_{\hypo_3}&=(\kappa_2(k))^{\sharp},\quad
(\kappa_1(k))^{\sharp}&=[2]_{\hypo_3}&= \kappa_2(k^{\sharp}) &\quad \text{if}\ k = i_3,\\
\kappa_1(k^{\sharp})&=[21]_{\hypo_3}&=(\kappa_2(k))^{\sharp}, \quad
(\kappa_1(k))^{\sharp}&=[\varepsilon]_{\hypo_3}&=\kappa_2(k^{\sharp}) & \quad\text{if}\  i_3< k< i_4,\\
\kappa_1(k^{\sharp})&=[1]_{\hypo_3}&=(\kappa_2(k))^{\sharp}, \quad
(\kappa_1(k))^{\sharp}&=[\varepsilon]_{\hypo_3}&=\kappa_2(k^{\sharp})& \quad \text{if}\ k=i_4,\\
\kappa_1(k^{\sharp})&=[\varepsilon]_{\hypo_3}&=(\kappa_2(k))^{\sharp},\quad
(\kappa_1(k))^{\sharp}&=[\varepsilon]_{\hypo_3}&=\kappa_2(k^{\sharp}) &\quad \text{otherwise}.
\end{array}\right.
\end{align*}}
Therefore it is routine to verify that $\kappa_{ij}(\bfw^{\sharp})=(\kappa_{ij}(\bfw))^{\sharp}$ for any $\bfw\in \mathcal{A}_n^{\star}$.

Now we can define the map $\varphi_{ij}: \mathcal{A}^{\star}_n \rightarrow \hypo_3\times \hypo_3$ by
\begin{align*}
\varphi_{ij}=
\left\{
\begin{array}{ll}
\lambda_{ij} & \hbox{if $i^{\sharp}=j$,} \\
[0.1cm]
\theta_{ij} & \hbox{if $i<j=j^{\sharp}<i^{\sharp}$ or $j^{\sharp}<i=i^{\sharp}<j$,} \\
[0.1cm]
\kappa_{ij} & \hbox{if $i<j^{\sharp}<j<i^{\sharp}$ or $j^{\sharp}<i<i^{\sharp}<j$}\\
\end{array}
\right.
\end{align*}
where $\lambda_{ij}, \theta_{ij}, \kappa_{ij}$ are defined as above. It follows from the previous analysis that $\varphi_{ij}$ is a homomorphism from $(\mathcal{A}_n^{\star},^{\sharp})$ to $(\hypo_3\times \hypo_3,^{\sharp})$.

\begin{lemma}\label{lemma:phi_factor_homomorphism}
The homomorphism $\varphi_{ij}$ induces a homomorphism $\varphi_{ij} : (\hypo_n,^{\sharp}) \rightarrow (\hypo_3\times \hypo_3,^{\sharp})$ for any $n\geq 4$.
\end{lemma}
	
\begin{proof}
To show that the homomorphism $\varphi_{ij}$ induces a homomorphism from $(\hypo_n,^\sharp)$ to $(\hypo_3\times \hypo_3,^{\sharp})$, we only need to show that for any $\bfu,\bfv\in \mathcal{A}_n^{\star}$, if $\bfu\equiv_{\hypo_\infty}\bfv$, then $\varphi_{ij}(\bfu)=\varphi_{ij}(\bfv)$.
It follows from Proposition \ref{pro:inversion} and the definition of $\varphi_{ij}$ that the evaluations of the first and the second components of $\varphi_{ij}(\bfu)$ and $\varphi_{ij}(\bfv)$ are the same respectively. In the following, we prove that the inversions of the first and the second components of $\varphi_{ij}(\bfu)$ and $\varphi_{ij}(\bfv)$ are the same respectively.

{\bf Case 1.} $\varphi_{ij}=\lambda_{ij}$. If there exists $k\in \con(\bfu)$ satisfying $i<k<j$, then it follows from $\ev(\bfu)=\ev(\bfv)$ that the first components of $\lambda_{ij}(\bfu)$ and $\lambda_{ij}(\bfv)$ have a $3$-$1$ inversion; if there is no $k\in\con(\bfu)$ satisfying $i<k<j$, then since the $3$-$1$ inversion in the first component of $\lambda_{ij}(\bfu)$ [resp. $\lambda_{ij}(\bfv)$] corresponds to the $j$-$i$ inversion in $\bfu$ [resp. $\bfv$], it follows from  Proposition \ref{pro:inversion} that the first component of $\lambda_{ij}(\bfu)$ has a $3$-$1$ inversion if and only if the first component of $\lambda_{ij}(\bfv)$ has a $3$-$1$ inversion. Therefore, the inversions of the first components of $\varphi_{ij}(\bfu)$ and $\varphi_{ij}(\bfv)$ are the same.  A similar argument can show that the second components of $\varphi_{ij}(\bfu)$ and $\varphi_{ij}(\bfv)$ are the same.

{\bf Case 2.} $\varphi_{ij}=\theta_{ij}$ or $\kappa_{ij}$.
Similar arguments with Case 1 can show that the first component of $\varphi_{ij}(\bfu)$ and $\varphi_{ij}(\bfv)$ have the same $2$-$1$ inversion and the second component of $\varphi_{ij}(\bfu)$ and $\varphi_{ij}(\bfv)$ have the same $3$-$2$ inversion.
\end{proof}
	
\begin{corollary}\label{coro:varphi_equiv_hypon}
Let $\bfu,\bfv \in \mathcal{A}_n^{\star}$ for any $n\geq 4$. Then $\bfu \equiv_{\hypo_\infty} \bfv$ if and only if $\varphi_{ij}(\bfu) = \varphi_{ij}(\bfv)$ for all $1 \leq i < j \leq n$.
\end{corollary}

\begin{proof}
The necessity follows from the proof of Lemma~\ref{lemma:phi_factor_homomorphism}. Let $\bfw \in \mathcal{A}_n^{\star}$ for some $n \geq 4$. Suppose $\con(\bfw)= \{a_1 < \dots < a_{\ell}\}$ for some $\ell \in \mathbb{N}$. For any $i$ with $1 \leq i < \ell$, if $a_i^{\sharp}=a_{i+1}$, then since there is no $k\in\con(\bfw)$ satisfying $a_i<k<a_{i+1}$, it follows from the definition of $\varphi_{ij}$ that the number of occurrences of $1$ in a component of $\varphi_{ij}(\bfw)$ equals to $\occ(a_i, \bfw)$ and the number of occurrences of $3$ in a component of $\varphi_{ij}(\bfw)$ equals to $\occ(a_{i+1}, \bfw)$ and the $3$-$1$ inversion in a component of $\varphi_{ij}(\bfw)$ corresponds to the $a_{i+1}$-$a_{i}$ inversion in $\bfw$; if
$a_i<a_{i+1}\leq a_{i+1}^{\sharp}<a_i^{\sharp}$ or $a_i<a_{i+1}^{\sharp}<a_{i+1}<a_i^{\sharp}$, then since there is no $k\in\con(\bfw)$ satisfying $a_i<k<a_{i+1}$, it follows from the definition of $\varphi_{ij}$ that the the number of occurrences of $1$ in the first component of $\varphi_{ij}(\bfw)$ equals to $\occ(a_i, \bfw)$ and the number of occurrences of $2$ in the first component of $\varphi_{ij}(\bfw)$ equals to $\occ(a_{i+1}, \bfw)$ and the $2$-$1$ inversion in the first component of $\varphi_{ij}(\bfw)$ corresponds to the $a_{i+1}$-$a_{i}$ inversion in $\bfw$; if $a_{i+1}^{\sharp}<a_{i}^{\sharp}\leq a_{i}<a_{i+1}$ or $a_{i+1}^{\sharp}<a_{i}<a_{i}^{\sharp}<a_{i+1}$, then since there is no $k\in\con(\bfw)$ satisfying $a_i<k<a_{i+1}$, it follows from the definition of $\varphi_{ij}$ that the number of occurrences of $2$ in the second component of $\varphi_{ij}(\bfw)$ equals to $\occ(a_i, \bfw)$ and the number of occurrences of $3$ in the second component of $\varphi_{ij}(\bfw)$ equals to $\occ(a_{i+1}, \bfw)$ and the $3$-$2$ inversion in the second component of $\varphi_{ij}(\bfw)$ corresponds to the $a_{i+1}$-$a_{i}$ inversion in $\bfw$.
Therefore both the number of occurrences of $a_i$ and $a_{i+1}$ in $\bfw$ and whether $\bfw$ has an $a_{i+1}$-$a_i$ inversion can be derived from the maps $\varphi_{a_i,a_{i+1}}$. Thus the sufficiency follows from Proposition~\ref{pro:inversion}.
\end{proof}

For each $n \in \mathbb{N}$, with $n \geq 4$, let $I_n$ be the index set
\[
\{ (i,j): 1 \leq i < j \leq n \}.
\]
Now, consider the map
\[
\phi_n : (\hypo_n,^{\sharp}) \rightarrow \prod\limits_{I_n} (\hypo_3\times \hypo_3,^{\sharp}),
\]
whose $(i,j)$-th component is given by $\varphi_{ij}([\bfw]_{\hypo_n})$ for $\bfw \in \mathcal{A}_n^{\star}$ and $(i,j) \in I_n$.
		
\begin{proposition}\label{prop:hypon_embed_hypo3}
The map $\phi_n$ is an embedding from $(\hypo_n,^{\sharp})$ to $\prod\limits_{I_n} (\hypo_3\times \hypo_3,^{\sharp})$.
\end{proposition}
	
\begin{proof}
The map $\phi_n$ is a homomorphism by Lemma~\ref{lemma:phi_factor_homomorphism}. It follows from the definition of $\phi_n$ and Corollary~\ref{coro:varphi_equiv_hypon} that $[\bfu]_{\hypo_n} =[\bfv]_{\hypo_n}$ if and only if $\phi_n([\bfu]_{\hypo_n}) = \phi_n([\bfv]_{\hypo_n})$ for any $\bfu,\bfv \in \mathcal{A}_n^{\star}$. Hence $\phi_n$ is an embedding.
\end{proof}

For any $([\bfu_1]_{\hypo_3}, [\bfu_2]_{\hypo_3}, \dots, [\bfu_{2|I_n|}]_{\hypo_3})\in \prod\limits_{2I_n}\hypo_3$, define an involution operation $^{\sharp}$ on $\prod\limits_{2I_n}\hypo_3$ by
\begin{align*}
([\bfu_1]_{\hypo_3}, [\bfu_2]_{\hypo_3}, \dots, [\bfu_{2|I_n|}]_{\hypo_3})^{\sharp}=([\bfu_{2|I_n|}]_{\hypo_3}^{\sharp}, \dots, [\bfu_{2}]_{\hypo_3}^{\sharp}, [\bfu_1]_{\hypo_3}^{\sharp}).
\end{align*}
Define a map $\eta:\prod\limits_{I_n} (\hypo_3\times \hypo_3,^{\sharp})\rightarrow (\prod\limits_{2I_n}\hypo_3,^{\sharp})$ given by
\begin{align*}
(([\bfu_1]_{\hypo_3}, [\bfu_2]_{\hypo_3}), ([\bfu_3]_{\hypo_3}, [\bfu_4]_{\hypo_3}), \dots, ([\bfu_{2|I_n|-1}]_{\hypo_3}, [\bfu_{2|I_n|}]_{\hypo_3}))\\
\mapsto ([\bfu_1]_{\hypo_3}, [\bfu_3]_{\hypo_3}, \dots, [\bfu_{2|I_n|-1}]_{\hypo_3}, [\bfu_{2|I_n|}]_{\hypo_3}, \dots, [\bfu_4]_{\hypo_3}, [\bfu_2]_{\hypo_3}).
\end{align*}
It is routine to verify that the map $\eta$ is an isomorphism. By Theorem \ref{thm:hypo3-repre}, each element in $(\hypo_3,^{\sharp})$ corresponds to a matrix in $UT_{13}(\mathbb{S})$ and the involution $^{\sharp}$ on $(\hypo_3,^{\sharp})$ corresponds to the skew transposition on $UT_{13}(\mathbb{S})$. It follows that there is an embedding, denoted by $\rho$, from $(\prod\limits_{2I_n}\hypo_3,^{\sharp})$ to $(UT_{26|I_n|}(\mathbb{S}),^{D})$. Let
\[
\psi_n=\rho\circ\eta\circ\phi_n.
\]
Then the following result holds.

\begin{theorem}\label{thm_hypo_tropical}
For each $n\geq 4$,  the map $\psi_n:(\hypo_n,^\sharp)\rightarrow(UT_{26|I_n|}(\mathbb{S}),^{D})$ is a faithful representation of $(\hypo_n,^\sharp)$.
\end{theorem}

\begin{remark}
Cain et al. have shown that $\hypo_n$ can be embedded into a direct product of $n$ copies of the free monogenic monoid and $\frac{n(n-1)}{2}$ copies of the finite monoid $A_0^1$ (denoted by $\mathcal{H}$ in \cite{CJKM21}) which will be defined in Section \ref{sec:characterization}. They gave a faithful representation of the free monogenic monoid as a monoid of $1\times 1$ matrices and a faithful representation of $A_0^1$ as a monoid of $2\times 2$ matrices. It follows that there is a faithful representation from $\hypo_n$ to $M_{n^2}(\mathbb{S})$ \cite[Theorems~3.3 and 3.4]{CJKM21}. Since the representation of $A_0^1$ they gave can not be extended to a representation from $(A_0^1, ^*)$ to $(M_{2}(\mathbb{S}),^{D})$,  the representation can not be extended to a representation from $(\hypo_n,^{\sharp})$ to $(M_{n^2}(\mathbb{S}),^{D})$.
\end{remark}

Let $a<b<c<d$ be a 4-element ordered alphabet and
\[
H=\langle a,b,c,d\,|\,\mathcal{R}_{\hypo_\infty}, ac=ca, ad=da, bc=cb, bd=db\rangle\cup \{1\}
\]
be a monoid.
The involution operation $^{\sharp}$ on $H$ can be defined by $a\mapsto d, b\mapsto c$.
By the definition of $H$ and Proposition \ref{pro:inversion}, the elements of $H$ are
characterized by their evaluation, their $b$-$a$ inversions and
their $d$-$c$ inversions.

\begin{theorem}\label{thm:same-var}
For any $m,n\geq 4$, the involution monoids $(\hypo_m,^{\sharp})$ and $(\hypo_n,^{\sharp})$ generate the same variety.
\end{theorem}
\begin{proof}
It follows from Proposition \ref{prop:hypon_embed_hypo3} that
\[
\var(\hypo_4,^{\sharp})\subseteq \var(\hypo_5,^{\sharp}) \subseteq \cdots \subseteq\var (\hypo_3\times \hypo_3,^{\sharp}).
\]
It suffices to show that $\var (\hypo_3\times \hypo_3,^{\sharp})\subseteq \var(\hypo_4,^{\sharp})$.
It is easy to see that $(H, ^{\sharp}$) is a homomorphism image of $(\hypo_4,^{\sharp})$, and so $\var(H,^{\sharp})\subseteq \var(\hypo_4,^{\sharp})$. In the following, we show that $\var (\hypo_3\times \hypo_3,^{\sharp})\subseteq \var(H,^{\sharp})$.

Let $\tau$ be a map from $(\hypo_3\times \hypo_3,^{\sharp})$ to $(H,^{\sharp})\times (H,^{\sharp})\times (H,^{\sharp})$ given by
\begin{align*}
&([\varepsilon]_{\hypo_3}, [\varepsilon]_{\hypo_3})\mapsto (1,1,1),\\
([1]_{\hypo_3},[\varepsilon]_{\hypo_3})&\mapsto (a,a,1), \quad \,\,([2]_{\hypo_3},[\varepsilon]_{\hypo_3})\mapsto (b, ba,a), \\ ([3]_{\hypo_3},[\varepsilon]_{\hypo_3})&\mapsto (1,b,b), \quad\, \,\,([\varepsilon]_{\hypo_3}, [1]_{\hypo_3})\mapsto (1,c,c), \\
([\varepsilon]_{\hypo_3}, [2]_{\hypo_3})&\mapsto (c,dc,d), \quad ([\varepsilon]_{\hypo_3}, [3]_{\hypo_3})\mapsto (d,d,1).
\end{align*}
We show that if $\bfu_1\equiv_{\hypo_{\infty}}\bfv_1, \bfu_2\equiv_{\hypo_{\infty}}\bfv_2$, then
\[
\tau([\bfu_1]_{\hypo_3},[\bfu_2]_{\hypo_3})=(U_1,U_2,U_3)=(V_1,V_2,V_3)=\tau([\bfv_1]_{\hypo_3},[\bfv_2]_{\hypo_3}).
\]
Since $\bfu_1\equiv_{\hypo_{\infty}}\bfv_1, \bfu_2\equiv_{\hypo_{\infty}}\bfv_2$, it is easy to see that $\ev(U_i)=\ev(V_i)$ for $i=1,2,3$.
It follows from Proposition \ref{pro:inversion} and the definition of $(H, ^{\sharp})$ that, for any $([\bfw_1]_{\hypo_3}, [\bfw_2]_{\hypo_3})\in\hypo_3\times \hypo_3$, the inversions of each component of $\tau([\bfw_1]_{\hypo_3}, [\bfw_2]_{\hypo_3})=(W_1, W_2, W_3)$  can be characterized as follows:
\begin{align*}
b\textrm{-}a \in\inver(W_1)~[\text{resp. } d\textrm{-}c \in\inver(W_3)]&\Leftrightarrow  2\textrm{-}1\in \inver(\bfw_1)~[\text{resp. } \inver(\bfw_2)],\\
b\textrm{-}a \in\inver(W_3)~[\text{resp. } d\textrm{-}c \in\inver(W_1)]&\Leftrightarrow 3\textrm{-}2\in \inver(\bfw_1)~[\text{resp. } \inver(\bfw_2)],\\
b\textrm{-}a \in\inver(W_2)~[\text{resp. } d\textrm{-}c \in\inver(W_2)] &\Leftrightarrow 3\textrm{-}1\in \inver(\bfw_1)~[\text{resp. } \inver(\bfw_2)].
\end{align*}
Since $\bfu_1\equiv_{\hypo_{\infty}}\bfv_1, \bfu_2\equiv_{\hypo_{\infty}}\bfv_2$, it is easy to see that $\inver(U_i)=\inver(V_i)$ for $i=1,2,3$. It is routine to check that the images of the generators of $(\hypo_3\times \hypo_3,^{\sharp})$
under $\tau$ followed by the involution of $(H,^{\sharp})$ are the same as those under the
involution of $(\hypo_3\times \hypo_3,^{\sharp})$ followed by $\tau$. Thus $\tau$ is a homomorphism.
By the previous observations and Proposition \ref{pro:inversion}, we can also conclude that $\tau$ is injective.
\end{proof}

\section{The identities satisfied by $(\hypo_n,^\sharp)$} \label{sec:characterization} %

In this section, we obtain a complete characterization of the word identities satisfied by $(\hypo_n,^\sharp)$ for each finite $n$.
Clearly, a word identity $\bfu\approx \bfv$ holds in $(\hypo_1,^\sharp)$ if and only if $\occ(x,\overline{\bfu})=\occ(x,\overline{\bfv})$ for any $x\in \con(\overline{\bfu\bfv})$.

Let $\bfu$ be a word. Define
\begin{align*}
\mix(\bfu)&=\{x\,|\, x, x^* \in \con(\bfu)\},\\
\ml(\bfu)\ &=\{x\,|\,x\in \mix(\bfu), \occ(x,\bfu)=1\},\\
\lin(\bfu)\ &=\{x \,|\, x\not\in\mix(\bfu), \occ(x,\bfu)=1\}.
\end{align*}
Clearly, $x\in \mix(\bfu)$ if and only if $x^*\in \mix(\bfu)$, $\mix(\bfu)\cap\lin(\bfu)= \emptyset$. For example, if $\bfu=x^*z^2xy^*x$, then $\mix(\bfu)=\{x, x^*\}, \ml(\bfu)=\{x^*\}$ and $\lin(\bfu)=\{y^*\}$.

Let
\[
A_0^1=\langle\, a, b\, |\, a^2 = a, b^2= b, aba =bab=ba\,\rangle\cup \{1\}=\{a, b, ab, ba, 1\}.
\]
The involution monoid $(A_0^1, ^*)$ can be defined by $A_0^1$ under the unary
operation $^*$ that interchanges the generators $a$ and $b$ and fixes all other elements. The involution monoid $(A_0^1, ^*)$ has been studied in \cite{GZL-A01} which showed that $(A_0^1, ^*)$ is the first example of a non-finitely based involution semigroup of order five. Clearly the generators of $A_0^1$ satisfy the relations in $\mathcal{R}_{\hypo_\infty}$, thus $(A_0^1, ^*)$ is a homomorphic image of $(\hypo_2,^\sharp)$.

\begin{theorem}\label{thm:A01}
A word identity $\bfu\approx \bfv$ holds in $(A_0^1, ^*)$ if and only if $\bfu\approx \bfv$ satisfies the following conditions:
\begin{enumerate}[\rm(i)]
\item $\con(\bfu)=\con(\bfv), \mix(\bfu)=\mix(\bfv), \lin(\bfu)=\lin(\bfv)$;
\item for any $x, y \in \con(\bfu)$,
\begin{enumerate}[\rm(a)]
  \item if $x,y\in \mix(\bfu)$, then $\{x, y\}\prec_{\bfu}\{x^*, y^*\}$ if and only if $\{x, y\}\prec_{\bfv}\{x^*, y^*\}$;
  \item if $x\in \mix(\bfu), y\not\in \mix(\bfu)$, then $x\prec_{\bfu}\{x^*, y\}$ if and only if $x \prec_{\bfv}\{x^*, y\}$ and $\{x, y\}\prec_{\bfu} x^* $ if and only if $\{x, y\} \prec_{\bfv} x^*$;
  \item if $x,y\not\in \mix(\bfu)$, then $x\prec_{\bfu} y$ if and only if $x\prec_{\bfv}y$.
\end{enumerate}
\end{enumerate}
\end{theorem}

\begin{proof}
Let $\bfu\approx \bfv$ be a word identity satisfied by $(A_0^1,^*)$. Then it is easy to show or see \cite[Lemma~4.2]{GZL-A01} that $\con(\bfu) = \con(\bfv)$. Suppose that $x \in \mix(\bfu)$. Then $x^* \in \mix(\bfu)$, and so $\{x, x^*\}\subseteq \con(\bfu) = \con(\bfv)$. This implies that $x \in \mix(\bfv)$ and so $\mix(\bfu) \subseteq \mix(\bfv)$. That $\mix(\bfv) \subseteq \mix(\bfu)$ by symmetry. Hence $\mix(\bfv) = \mix(\bfu)$.
Suppose that there exists $x\in \lin(\bfu) \setminus \lin(\bfv)$. Then $x\not \in \mix(\bfu)=\mix(\bfv)$ and $\occ(x,\bfv)\geq 2$. Let $\varphi_1$ be the homomorphism that maps $x$ to $ab$ and any other variable to $1$. Then $\varphi_1(\bfu)=ab\neq ba=\varphi_1(\bfv)$, a contradiction. Hence $\lin(\bfu)=\lin(\bfv)$. Therefore the condition (i) holds.

Suppose that there exist $x, y \in \mix(\bfu)$ such that $\{x, y\}\prec_{\bfu}\{x^*, y^*\}$ but $\{x, y\}\not\prec_{\bfv}\{x^*, y^*\}$. Then $\bfu[x, y] \in \{x,y\}^{+} \cdot \{x^*,y^*\}^{+}$ while either some $x$ occurs after the first $x^*$ or the first $y^*$, or some $y$ occurs after the first $x^*$ or the first $y^*$ in $\bfv$. Let $\varphi_2$ be the homomorphism that maps $x, y$ to $a$ and any other variable to $1$. Then $\varphi_2(\bfu)=ab\neq ba=\varphi_2(\bfv)$, a contradiction. Hence $\{x, y\}\prec_{\bfu}\{x^*, y^*\}$ if and only if $\{x, y\}\prec_{\bfv}\{x^*, y^*\}$ for any $x, y\in \mix(\bfu)$. Thus (iia) holds.
Suppose that there exist $x\in \mix(\bfu), y\not\in \mix(\bfu)$ such that $\{x, y\}\prec_{\bfu}x^* $ but $\{x, y\}\not\prec_{\bfv}x^*$. Then $\bfu[x, y] \in \{x,y\}^{+} \cdot \{x^*\}^{+}$ while either some $x$ occurs after the first $x^*$, or some $y$ occurs after the first $x^*$ in $\bfv$. Thus $\varphi_2(\bfu)=ab\neq ba =\varphi_2(\bfv)$, a contradiction. Hence $\{x, y\}\prec_{\bfu} x^*$ if and only if $\{x, y\}\prec_{\bfv} x^*$ for any $x\in \mix(\bfu), y\not\in \mix(\bfu)$.
Suppose that there exist $x\in \mix(\bfu)$ and $y\not\in \mix(\bfu)$ such that $x\prec_{\bfu}\{x^*, y\}$ but $x\not\prec_{\bfv}\{x^*, y\}$. Then $\bfu[x, y] \in \{x\}^{+} \cdot \{x^*,y\}^{+}$ while some $x$ occurs after the first $x^*$ or the first $y$ in $\bfv$. Let $\varphi_3$ be the homomorphism that maps $x$ to $a$, $y$ to $b$ and any other variable to $1$. Then $\varphi_3(\bfu)=ab\neq  ba=\varphi_3(\bfv)$, a contradiction. Hence $x \prec_{\bfu}\{x^*, y^*\}$ if and only if $x \prec_{\bfv}\{x^*, y^*\}$ for any $x\in\mix(\bfu), y\not\in \mix(\bfu)$. Thus (iib) holds.
Suppose that there exist $x, y \not\in \mix(\bfu)$ such that $x\prec_{\bfu} y$ but $x\not\prec_{\bfv}y$. Then $\bfu[x, y] \in \{x\}^{+} \cdot \{y\}^{+}$ while some $x$ occurs after the first $y$ in $\bfv$. Thus $\varphi_3(\bfu)=ab\neq ba=\varphi_3(\bfv)$, a contradiction. Hence $x\prec_{\bfu} y$ if and only if $x\prec_{\bfv}y$ for any $x, y \not\in \mix(\bfu)$. Thus (iic) holds.
 Therefore the condition (ii) holds.

Conversely, let $\bfu\approx \bfv$ ba a word identity satisfying conditions (i) and (ii) and $\phi$ be any homomorphism from $(\mathcal{X}\cup \mathcal{X}^*)^{+}$ to $(A_0^1,^*)$.
First we show that $\phi(\bfu)=ab$ if and only if $\phi(\bfv)=ab$.
By symmetry, we may assume that $\phi(\bfu)=ab$.
Note that $\phi$ mapping a variable $x$ in $\bfu\approx \bfv$ to $1$ is the same as removing all occurrences of $x$ and $x^*$ in $\bfu\approx \bfv$.
Hence if $\phi(x)=1$ for some $x\in \con(\bfu)$, then we only need to consider the identity obtained from $\bfu\approx \bfv$ by deleting all occurrences of $x$ and $x^*$.
Therefore we may assume that $\phi(x)\neq 1$ for any $x\in \con(\bfu)$.
Now by the definition of $A_0^1$, either $\bfu=\bfu_1\bfu_2$ for some nonempty words $\bfu_1, \bfu_2$ such that $\phi(\bfu_1)=a, \phi(\bfu_2)=b$, or $\bfu=\bfu_1z\bfu_2$ for some possibly empty words $\bfu_1, \bfu_2$ such that $\phi(\bfu_1)=a$, $\phi(\bfu_2)=b$ and $\phi(z)=ab$.
If $\bfu=\bfu_1\bfu_2$, then since $\phi(x)=a$ for any $x\in \mathsf{con}(\bfu_1)$ and $\phi(y)=b$ for any $y\in \mathsf{con}(\bfu_2)$, it follows that $\con(\bfu_1)\cap \con(\bfu_2) =\emptyset$.
Therefore $\con(\bfu_1) \prec_{\bfu}\con(\bfu_2)$.
Suppose that $ \con(\bfu_1) \not\prec_{\bfv} \con(\bfu_2)$.
Then there exist some $x\in \con(\bfu_1)$ and $y\in \con(\bfu_2)$ such that $_1y\prec_{\bfv} {_\infty x}$.
If $x,y \in \mix(\bfu)$, then $x^* \in \con(\bfu_2)$ and $y^* \in \con(\bfu_1)$ satisfying $\{x,y^{*}\} \prec_{\bfu} \{x^*,y\}$, whence $\{x,y^{*}\} \prec_{\bfv} \{x^*,y\}$ by (iia);
if $x \in \mix(\bfu)$ and $y \notin \mix(\bfu)$, then $x^* \in \con(\bfu_2)$ satisfying $x \prec_{\bfu} \{x^*,y\}$, whence $x \prec_{\bfv} \{x^*,y\}$ by (iib);
if $x \notin \mix(\bfu)$ and $y \in \mix(\bfu)$, then $y^* \in \con(\bfu_1)$ satisfying $\{x,y^*\} \prec_{\bfu} y$, whence $\{x, y^*\} \prec_{\bfv} y$ by (iib);
if $x \notin \mix(\bfu)$ and $y \notin \mix(\bfu)$, then $x \prec_{\bfu} y$, whence $x \prec_{\bfv} y$ by (iic). Hence in any case, ${_\infty x} \prec_{\bfv} {_1y}$, a contradiction.
Therefore $ \con(\bfu_1) \prec_{\bfv} \con(\bfu_2)$.
Consequently, $\phi(\bfv)=ab$. If $\bfu=\bfu_1z\bfu_2$, then $\phi(x)=a$ for any $x\in \mathsf{con}(\bfu_1)$, $\phi(y)=b$ for any $y\in \mathsf{con}(\bfu_2)$, and $z, z^* \not \in \con(\bfu_1\bfu_2)$.
Hence $\con(\bfu_1) \cap \con(\bfu_2) =\emptyset$ and $z\in \lin(\bfu)$. Therefore, $\con(\bfu_1)\prec_{\bfu} z \prec_{\bfu}\con(\bfu_2)$.
Suppose that $\con(\bfu_1)\not\prec_{\bfv} z$. Then there exists $x\in \con(\bfu_1)$ such that $_1 z\prec_{\bfv} {_\infty x}$.
Since $x \prec_{\bfu} z$, it follows from (iib) and (iic) that $x \prec_{\bfv} z$, a contradiction.
Hence $\con(\bfu_1)\prec_{\bfv} z$. Suppose that $z\not\prec_{\bfv} \con(\bfu_2)$.
Then there exists $x\in \con(\bfu_2)$ such that $_1x \prec_{\bfv} {_\infty z}$.
Since $z \prec_{\bfu} x$, it follows from (iib) and (iic) that $z \prec_{\bfv} x$, a contradiction.
Hence $ z \prec_{\bfv} \con(\bfu_2)$. Therefore $\con(\bfu_1)\prec_{\bfv} z \prec_{\bfv}\con(\bfu_2)$. Consequently $\phi(\bfv)=ab$.

Clearly $\phi(\bfu)=1$ [resp. $a, b$] if and only if $\phi(\bfv)=1$ [resp. $a, b$] by the definition of $A_0^1$ and the condition (i). These imply that $\phi(\bfu)=ba$ if and only if $\phi(\bfv)=ba$. Therefore any word identity $\bfu\approx \bfv$ satisfying conditions (i) and (ii) holds in $(A_0^1,^*)$.
\end{proof}

Let
\[
A=\langle\, a, b\, |\, ab=ba\,\rangle=\{a^mb^n\,|\,m,n\geq 0\}
\]
 be a monoid.
The monoid $A$ forms an involution monoid $(A, ^*)$ under the unary operation $^*: a^mb^n\mapsto a^nb^m$.   It is routine to verify that $(A,^*)$ is a homomorphic image of $(\hypo_2,^\sharp)$ under the map given by $[1]_{\hypo_2}\mapsto a$ and $[2]_{\hypo_2}\mapsto b$. A word identity $\bfu \approx \bfv$ is \textit{balanced} if $\mathsf{occ}(x, \bfu)=\mathsf{occ}(x, \bfv)$ for any $x\in \mathsf{con}(\bfu\bfv)$.

\begin{theorem}\label{thm:balanced}
A word identity $\bfu\approx \bfv$ holds in $(A, ^*)$ if and only if $\bfu\approx \bfv$ is balanced.
\end{theorem}

\begin{proof}
Suppose that $\bfu\approx \bfv$ is a word identity satisfied by $(A, ^*)$ such that $\occ(x,\bfu)\ne\occ(x,\bfv)$ or $\occ(x^*,\bfu)\ne \occ(x^*,\bfv)$ for some $x\in \mathcal{X}$. Let $\varphi$ be a homomorphism from $(\mathcal{X}\cup \mathcal{X}^*)^{+}$ to $(A,^*)$ that maps $x$ to $a$ and any other variable to $1$.
Then $\varphi(\bfu)=a^{\occ(x,\bfu)}b^{\occ(x^*,\bfu)} \ne a^{\occ(x,\bfv)}b^{\occ(x^*,\bfv)}=\varphi(\bfv)$, a contradiction.

Conversely, if $\bfu\approx \bfv$ is balanced, then $\varphi(\bfu)=\varphi(\bfv)$ for any homomorphism $\varphi$ from $(\mathcal{X}\cup \mathcal{X}^*)^{+}$ to $(A,^*)$ since $(A,^*)$ is commutative. Therefore any balanced word identity $\bfu\approx \bfv$ holds in $(A, ^*)$.
\end{proof}

\begin{theorem}\label{thm:hypo2}
A word identity $\bfu\approx \bfv$ holds in $(\hypo_2, ^{\sharp})$ if and only if
\begin{enumerate}[\rm(i)]
\item $\bfu\approx \bfv$ is balanced;
\item for any $x, y \in \con(\bfu)$,
\begin{enumerate}[\rm(a)]
  \item if $x,y\in \mix(\bfu)$, then $\{x, y\}\prec_{\bfu}\{x^*, y^*\}$ if and only if $\{x, y\}\prec_{\bfv}\{x^*, y^*\}$;
  \item if $x\in \mix(\bfu), y\not\in \mix(\bfu)$, then $x\prec_{\bfu}\{x^*, y\}$ if and only if $x \prec_{\bfv}\{x^*, y\}$ and $\{x, y\}\prec_{\bfu} x^* $ if and only if $\{x, y\} \prec_{\bfv} x^*$;
  \item if $x,y\not\in \mix(\bfu)$, then $x\prec_{\bfu} y$ if and only if $x\prec_{\bfv}y$.
\end{enumerate}
\end{enumerate}
\end{theorem}

\begin{proof}
Note that $(A,^*), (A_0^1,^*)\in\var(\hypo_2, ^{\sharp})$. Then $(A,^*)\times (A_0^1,^*) \in \var(\hypo_2, ^{\sharp})$. On the other hand,
it follows from Theorem~\ref{thm:hypo2-repre} that $(\hypo_2, ^{\sharp})$ is isomorphic to the submonoid of $(UT_5(\mathbb{S}), ^{D})$ generated by $\mathrm{E}_5, \mathsf{diag}\{s,\mathrm{J}, \mathbf{1}\}$ and $\mathsf{diag}\{\mathbf{1},\mathrm{K}, s\}$.
Note that $(A, ^*)$ is isomorphic to the involution matrix monoid generated by $\mathrm{P}, \mathrm{Q}, \mathrm{E}_2$ under the skew transposition and $(A_0^1, ^*)$ is isomorphic to the involution matrix monoid generated by $\mathrm{J}, \mathrm{K}, \mathrm{E}_3$ under the skew transposition.
Let $\varphi$ be a map from the involution matrix monoid generated by $\mathrm{E}_5, \mathsf{diag}\{s,\mathrm{J}, \mathbf{1}\}$ and $\mathsf{diag}\{\mathbf{1},\mathrm{K}, s\}$ to $(A,^*)\times (A_0^1,^*)$ given by
\[
\mathrm{E}_5\mapsto (\mathrm{E}_2,\mathrm{E}_3), \mathsf{diag}\{s,\mathrm{J}, \mathbf{1}\} \mapsto (\mathrm{P}, \mathrm{J}),  \mathsf{diag}\{\mathbf{1},\mathrm{K}, s\}\mapsto(\mathrm{Q}, \mathrm{K}).
\]
Then $\varphi$ is an embedding. Hence $(\hypo_2, ^{\sharp})\in \var((A,^*)\times (A_0^1,^*))$. Therefore $\var (\hypo_2, ^{\sharp})=\var((A,^*)\times (A_0^1,^*))$. Now the result follows from Theorems~\ref{thm:A01} and \ref{thm:balanced}.
\end{proof}

Let
\begin{align*}
B&=\left\langle\, a, b, c\,\, \begin{tabular}{|@{\,\,}c} $a^2 = a, b^2 = b, c^2=c, aba=bab=ba, aca=cac=ca,$\\[.1ex] $bcb=cbc=cb, acb=cab, bac=bca, bcab=cba$\end{tabular} \right\rangle\cup \{1\}\\
&=\{1, a,b,c, ab,ba, ac,ca, bc,cb, abc, acb, bac, cba\}.
\end{align*}
The monoid $B$ admits a unique involution $^*$ which is induced by the mapping $(a, b, c) \mapsto (c,b,a)$.  Clearly the generators of $B$ satisfy the relations in $\mathcal{R}_{\hypo_\infty}$, thus $(B, ^*)$ is a homomorphic image of $(\hypo_3,^\sharp)$.

\begin{theorem}\label{thm:B}
A word identity $\bfu\approx \bfv$ holds in $(B, ^*)$ if and only if $\bfu\approx \bfv$ satisfies the conditions
\begin{enumerate}[\rm(i)]
\item $\con(\bfu)=\con(\bfv), \mix(\bfu)=\mix(\bfv), \lin(\bfu)=\lin(\bfv)$;
\item for any $x, y \in \con(\bfu)$,
\begin{enumerate}[\rm(a)]
  \item if $x,y\in \mix(\bfu)$, then $\{x, y\}\prec_{\bfu}\{x^*, y^*\}$ if and only if $\{x, y\}\prec_{\bfv}\{x^*, y^*\}$;
  \item if $x\in \mix(\bfu), y\not\in \mix(\bfu)$, then $x\prec_{\bfu}\{x^*, y\}$ if and only if $x \prec_{\bfv}\{x^*, y\}$ and $\{x, y\}\prec_{\bfu} x^* $ if and only if $\{x, y\} \prec_{\bfv} x^*$;
  \item if $x,y\not\in \mix(\bfu)$, then $x\prec_{\bfu} y$ if and only if $x\prec_{\bfv}y;$
\end{enumerate}

\item if $x\in \mix(\bfu)$ and $x\prec_{\bfu} x^*$, then $x \in \ml(\bfu)$ if and only if $x \in \ml(\bfv)$ and $x^* \in \ml(\bfu)$ if and only if $x^* \in \ml(\bfv)$;
\item if $x\in \mix(\bfu)$ and $y \in \con(\bfu)$, then $\{x, x^*\}\prec_{\bfu} y$ if and only if $\{x, x^*\}\prec_{\bfv} y$ and $y \prec_{\bfu}\{x, x^*\}$ if and only if $y \prec_{\bfv}\{x, x^*\}$.
\end{enumerate}
\end{theorem}

\begin{proof}
Let $\bfu\approx \bfv$ be a word identity satisfied by $(B, ^*)$.
Clearly, the involution submonoid of $(B, ^*)$ generated by $a, c$ is isomorphic to $(A_0^1, ^*)$.  It follows from Theorem \ref{thm:A01} that conditions (i) and (ii) hold. If $x \in \mix(\bfu)$ and $x\prec_{\bfu} x^*$, then $\bfu[x] \in x^{\alpha}(x^*)^{\beta}$ for some $\alpha, \beta \geq 1$. It follows from the condition (iia) that $\bfv[x] \in x^{\alpha'}(x^*)^{\beta'}$ for some $\alpha', \beta' \geq 1$. Suppose that $\alpha=1$ but $\alpha'>1$. Let $\varphi_1$ be the homomorphism from $(\mathcal{X}\cup \mathcal{X}^*)^{+}$ to $(B,^*)$ that maps $x$ to $ab$ and any other variable to $1$. Then $\varphi_1(\bfu)\in \{abc, acb\}, \varphi_1(\bfv)\in \{bac, cba\}$, which implies $\varphi_1(\bfu) \ne \varphi_1(\bfv)$, a contradiction. Hence $\alpha=1$ if and only if $\alpha'=1$. A similar argument can show that $\beta=1$ if and only if $\beta'=1$. Therefore the condition (iii) holds.

Suppose that there exist $x\in \mix(\bfu)$ and $y \in\con(\bfu)$ such that $\{x, x^*\}\prec_{\bfu}  y$ but $\{x, x^*\}\not\prec_{\bfv} y$. Then $\bfu[x, y] \in \{x,x^*,y^*\}^{+}\{y,y^*\}^{+}$ and some $x$ or $x^*$ occurs after the first $y$ in $\bfv$.  Let $\varphi_2$ be the homomorphism from $(\mathcal{X}\cup \mathcal{X}^*)^{+}$ to $(B, ^*)$ that maps $x$ to $b$, $y$ to $c$ and any other variable to $1$. Then
$\varphi_2(\bfu)\in \{bc, abc, bac\}$ and $\varphi_2(\bfv)\in \{cb, acb, cba\}$, which implies that $\varphi_2(\bfu) \ne \varphi_2(\bfv)$, a contradiction. Hence $\{x, x^*\}\prec_{\bfu} y$ if and only if $\{x, x^*\}\prec_{\bfv}y$. A similar argument can show that $y \prec_{\bfu}\{x, x^*\}$ if and only if $y \prec_{\bfv}\{x, x^*\}$. Therefore the condition (iv) holds.

Conversely, let $\bfu\approx \bfv$ be an identity that satisfies the conditions (i)--(iv) and $\phi$ be any homomorphism from $(\mathcal{X}\cup \mathcal{X}^*)^{+}$ to $(B,^*)$.
Since $(A_0^1, ^*)$ is isomorphic to the involution subsemigroup $\{a, c, ac, ca, 1\}$ of $(B, ^*)$, it follows from Theorem~\ref{thm:A01} that $\bfu\approx \bfv$ holds in $\{a, c, ac, ca, 1\}$. Hence $\phi(\bfu) = \phi(\bfv)$ when $\phi(\bfu) \in \{a, c, ac, ca, 1\}$.
Clearly $\phi(\bfu)=b$ if and only if $\phi(\bfv)=b$ by the definition of $B$ and the condition (i).

Next we have to show that $\phi(\bfu)=ab$ if and only if $\phi(\bfv)=ab$. By symmetry, we may assume that $\phi(\bfu)=ab$. Note that $\phi$ mapping a variable $x$ in $\bfu\approx \bfv$ to $1$ is the same as removing all occurrences of $x$ and $x^*$ in $\bfu\approx \bfv$.
Hence if $\phi(x)=1$ for some $x\in \con(\bfu)$, then we only need to consider the identity obtained from $\bfu\approx \bfv$ by deleting all occurrences of $x$ and $x^*$.
Therefore we may assume that $\phi(x)\neq 1$ for any $x\in \con(\bfu)$. Now by the definition of the monoid $B$, either $\bfu=\bfu_1\bfu_2$ for some nonempty words $\bfu_1, \bfu_2$ such that $\phi(\bfu_1)=a, \phi(\bfu_2)=b$, or $\bfu=\bfu_1z\bfu_2$ for some possibly empty words $\bfu_1, \bfu_2$ such that $\phi(\bfu_1)=a$, $\phi(\bfu_2)=b$ and $\phi(z)=ab$.
If $\bfu=\bfu_1\bfu_2$, then $\phi(x)=a$ and $x^*\not\in \con(\bfu)$ for any $x\in \mathsf{con}(\bfu_1)$ and $\phi(y)=b$ for any $y\in \mathsf{con}(\bfu_2)$. Hence $\con(\bfu_1) \cap \mix(\bfu)=\emptyset$ and $\con(\bfu_1)\cap \con(\bfu_2) =\emptyset$, whence $\con(\bfu_1) \prec_{\bfu}\con(\bfu_2)$. Suppose that $ \con(\bfu_1) \not\prec_{\bfv} \con(\bfu_2)$.
Then there exist some $x\in \con(\bfu_1)$ and $y\in \con(\bfu_2)$ such that $_1y\prec_{\bfv} {_\infty x}$. If $y \in \mix(\bfu)$, then $y^* \in \con(\bfu_2)$ satisfying $x \prec_{\bfu} \{y, y^*\}$, whence $x \prec_{\bfv} \{y, y^*\}$ by (iv); if $y \notin \mix(\bfu)$, then  $x \prec_{\bfu} y$, whence $x \prec_{\bfv} y$ by (iic). Hence in any case, $_\infty x\prec_{\bfv} {_1y}$, a contradiction. Therefore $ \con(\bfu_1) \prec_{\bfv} \con(\bfu_2)$. It follows from (i) that $\phi(\bfv)=ab$. If $\bfu=\bfu_1z\bfu_2$, then $\phi(x)=a$ and $x^* \notin \con(\bfu)$ for any $x\in \mathsf{con}(\bfu_1)$, $\phi(y)=b$ for any $y\in \mathsf{con}(\bfu_2)$, and $z, z^* \not \in \con(\bfu_1\bfu_2)$. Hence $\con(\bfu_1) \cap \mix(\bfu)=\emptyset$, $\con(\bfu_1) \cap \con(\bfu_2) =\emptyset$ and $z\in \lin(\bfu)$, whence $\con(\bfu_1)\prec_{\bfu} z \prec_{\bfu}\con(\bfu_2)$. Suppose that $\con(\bfu_1)\not\prec_{\bfv} z$. Then there exists $x\in \con(\bfu_1)$ such that $_1z\prec_{\bfv} {_\infty x}$. Since $x \prec_{\bfu} z$, it follows from (iic) that $x \prec_{\bfv} z$, a contradiction. Hence $\con(\bfu_1)\prec_{\bfv} z$. Suppose that $z\not\prec_{\bfv} \con(\bfu_2)$. Then there exists $y\in \con(\bfu_2)$ such that $_1y \prec_{\bfv} {_\infty z}$. Since $z \prec_{\bfu} y$, it follows from (iic) and (iv) that $z \prec_{\bfv} y$, a contradiction. Hence $ z \prec_{\bfv} \con(\bfu_2)$. Therefore $\con(\bfu_1)\prec_{\bfv} z \prec_{\bfv}\con(\bfu_2)$. It follows from (i) that $\phi(\bfv)=ab$. A similar argument can show that $\phi(\bfu)=bc$ if and only if $\phi(\bfv)=bc$.

It remains to show that $\phi(\bfu)=ba$ if and only if $\phi(\bfv)=ba$. By symmetry, we may assume that $\phi(\bfu)=ba$. Then by the definition of $B$ and our assumption,  $\phi(x)\in \{a, b, ab, ba\}$ for any $x\in \con(\bfu)$. It follows from definition of $B$ and the condition (i) that $\phi(\bfv) \in \{a, b, ab, ba\}$. Note that $\phi(\bfu)=a$ [resp. $b, ab$] if and only if $\phi(\bfv)=a$ [resp. $b, ab$]. Hence $\phi(\bfv)=ba$. Similarly, $\phi(\bfu)=cb$ if and only if $\phi(\bfv)=cb$.

Now we show that $\phi(\bfu)=abc$ [resp. $acb, bac, cba$] if and only if $\phi(\bfv)=abc$ [resp. $acb, bac, cba$].
We only need to show that $a$ precedes $b$ in $\phi(\bfu)$ if and only if $a$ precedes $b$ in $\phi(\bfv)$ and that $b$ precedes $c$ in $\phi(\bfu)$ if and only if $b$ precedes $c$ in $\phi(\bfv)$.
In the following, we show that $a$ precedes $b$ in $\phi(\bfu)$ if and only if $a$ precedes $b$ in $\phi(\bfv)$, and $b$ precedes $c$ in $\phi(\bfu)$ if and only if $b$ precedes $c$ in $\phi(\bfv)$ can be obtained by a similar argument.
Without loss of generality we may assume that $a$ precedes $b$ in $\phi(\bfu)$.
Let $\mathcal{X}_1, \mathcal{X}_2, \mathcal{X}_3\subseteq \mathcal{X}$ satisfying $\phi(x)\in \{a, ac, ca\}$ for $x \in \mathcal{X}_1$, $\phi(y)\in \{b, bc, cb\}$ for $y \in \mathcal{X}_2$ and $\phi(w)=c$ for $w \in \mathcal{X}_3$.
It follows from the definition of $B$ that either $\bfu \in \{\mathcal{X}_1, \mathcal{X}_3\}^{+}\cdot \{\mathcal{X}_2, \mathcal{X}_3\}^{+}$, or $\bfu \in \{\mathcal{X}_1, \mathcal{X}_3\}^{+}\cdot z\cdot \{\mathcal{X}_2, \mathcal{X}_3\}^{+}$ satisfying $\phi(z)\in \{ab, abc, acb\}$.
If $\bfu \in \{\mathcal{X}_1, \mathcal{X}_3\}^{+}\cdot \{\mathcal{X}_2, \mathcal{X}_3\}^{+}$, then $\mathcal{X}_1\cap \mathcal{X}_2 =\emptyset$, whence $\mathcal{X}_1 \prec_{\bfu} \mathcal{X}_2$.
Suppose that $ \mathcal{X}_1 \not\prec_{\bfv} \mathcal{X}_2$. Then there exist some $x\in \mathcal{X}_1$ and $y\in \mathcal{X}_2$ such that $_1y\prec_{\bfv} {_\infty x}$.
Note that if $y \in \mix(\bfu)$, then $y^* \in \mathcal{X}_2$ satisfying $x \prec_{\bfu} \{y, y^*\}$, whence $x \prec_{\bfv} \{y, y^*\}$ by (iv);
if $y\notin \mix(\bfu)$, then  $x \prec_{\bfu} y$, whence $x \prec_{\bfv} y$ by (iic).
Hence in any case, $_\infty x \prec_{\bfv} {_1y}$, a contradiction.
Therefore $\mathcal{X}_1 \prec_{\bfv} \mathcal{X}_2$ and so $a$ precedes $b$ in $\phi(\bfv)$ by (i).
If $\bfu \in \{\mathcal{X}_1, \mathcal{X}_3\}^{+}\cdot z\cdot \{\mathcal{X}_2, \mathcal{X}_3\}^{+}$, then $z \not \in \mathcal{X}_1\mathcal{X}_2\mathcal{X}_3$, $z^*\in \mathcal{X}_2$ when $\phi(z)=ab$ and $z^*\not\in \mathcal{X}_1\mathcal{X}_2\mathcal{X}_3$ otherwise.
Hence $\mathcal{X}_1 \cap \mathcal{X}_2 =\emptyset$, whence $\mathcal{X}_1\prec_{\bfu} z \prec_{\bfu} \mathcal{X}_2$. Suppose that $\mathcal{X}_1\not\prec_{\bfv} z$.
Then there exists $x\in \mathcal{X}_1$ such that $_1z \prec_{\bfv} {_\infty x}$.
If $z\not\in \mix(\bfu)$, then $x \prec_{\bfv} z$ by $x \prec_{\bfu} z$ and (iic), a contradiction.
If $z\in \mix(\bfu)$, then $x \prec_{\bfv} \{z,z^*\}$ by $x \prec_{\bfu} \{z, z^*\}$ and (iv), a contradiction.
Hence $\mathcal{X}_1\prec_{\bfv} z$. Suppose that $z\not\prec_{\bfv} \mathcal{X}_2$.
Then there exists $y\in \mathcal{X}_2$ such that $_1y\prec_{\bfv} {_\infty z}$.
If $y=z^*$, then $z \prec_{\bfv} z^*$ by $z \prec_{\bfu} z^*$ and (iia), a contradiction.
If $y\neq z^*$ and $y\in \mix(\bfu)$, then $y^* \in \mathcal{X}_2$ satisfying $z \prec_{\bfu} \{y, y^*\}$, it follows from (iv) that $z \prec_{\bfv} \{y, y^*\}$, a contradiction.
If $y\neq z^*$ and $y\not \in \mix(\bfu)$, then $z \prec_{\bfv} y$ by $z \prec_{\bfu} y$ and (iic), a contradiction.
Hence $z \prec_{\bfv}\mathcal{X}_2$. Therefore $\mathcal{X}_1\prec_{\bfv} z \prec_{\bfv}\mathcal{X}_2$, and so $a$ precedes $b$ in $\phi(\bfv)$ by (i).

Therefore any word identity $\bfu\approx \bfv$ satisfying conditions (i)--(iv) is satisfied by $(B,^*)$.
\end{proof}

\begin{theorem}\label{hypo3}
A word identity $\bfu\approx \bfv$ holds in $(\hypo_3, ^{\sharp})$ if and only if
\begin{enumerate}[\rm(i)]
\item $\bfu\approx \bfv$ is balanced;
\item for any $x, y \in \con(\bfu)$,
\begin{enumerate}[\rm(a)]
  \item if $x,y\in \mix(\bfu)$, then $\{x, y\}\prec_{\bfu}\{x^*, y^*\}$ if and only if $\{x, y\}\prec_{\bfv}\{x^*, y^*\}$;
  \item if $x\in \mix(\bfu), y\not\in \mix(\bfu)$, then $x\prec_{\bfu}\{x^*, y\}$ if and only if $x \prec_{\bfv}\{x^*, y\}$ and $\{x, y\}\prec_{\bfu} x^* $ if and only if $\{x, y\} \prec_{\bfv} x^*$;
  \item if $x,y\not\in \mix(\bfu)$, then $x\prec_{\bfu} y$ if and only if $x\prec_{\bfv}y;$
\end{enumerate}

\item if $x\in \mix(\bfu)$ and $y \in \con(\bfu)$, then $\{x, x^*\}\prec_{\bfu} y$ if and only if $\{x, x^*\}\prec_{\bfv} y$ and $y \prec_{\bfu}\{x, x^*\}$ if and only if $y \prec_{\bfv}\{x, x^*\}$.
\end{enumerate}
\end{theorem}

\begin{proof}
Note that $(A,^*), (B,^*)\in\var(\hypo_3, ^{\sharp})$. Then $(A,^*)\times (B,^*)\in \var(\hypo_3, ^{\sharp})$. On the other hand,
it follows from Theorem~\ref{thm:hypo3-repre} that $(\hypo_3, ^{\sharp})$ is isomorphic to the submonoid of $(UT_{13}(\mathbb{S}), ^{D})$ generated by matrices $\mathrm{E}_{13}$, $\mathsf{diag}\{\mathrm{P},\mathrm{J},\mathrm{J},\mathrm{E}_3,\mathrm{E}_2\}$, $\mathsf{diag}\{\mathrm{Q},\mathrm{K},\mathrm{KJ},\mathrm{J},\mathrm{P}\}$ and $\mathsf{diag}\{\mathrm{E}_2,\mathrm{E}_3, \mathrm{K},\mathrm{K},\mathrm{Q}\}$. Note that $(A, ^*)$ is isomorphic to the involution matrix monoid generated by $\mathrm{P}, \mathrm{Q}, \mathrm{E}_2$ under the skew transposition and $(B, ^*)$ is isomorphic to the involution matrix monoid generated by $\mathrm{E}_{9}$, $\mathsf{diag}\{\mathrm{J},\mathrm{J},\mathrm{E}_3\}, \mathsf{diag}\{\mathrm{K},\mathrm{KJ},\mathrm{J}\}$, $\mathsf{diag}\{\mathrm{E}_3, \mathrm{K},\mathrm{K}\}$ under the skew transposition.
Let $\varphi$ be a map from the involution matrix monoid generated by $\mathrm{E}_{13}, \mathsf{diag}\{\mathrm{P},\mathrm{J},\mathrm{J},\mathrm{E}_3,\mathrm{E}_2\}$,  $\mathsf{diag}\{\mathrm{Q},\mathrm{K},\mathrm{KJ},\mathrm{J},\mathrm{P}\}$ and $\mathsf{diag}\{\mathrm{E}_2,\mathrm{E}_3, \mathrm{K}$, $\mathrm{K},\mathrm{Q}\}$ to $(A,^*)\times (B,^*)\times(A,^*)$ given by
\[
\begin{array}{lll}
&\mathrm{E}_{13}\mapsto (\mathrm{E}_2, 1, \mathrm{E}_2), & \mathsf{diag}\{\mathrm{P},\mathrm{J},\mathrm{J},\mathrm{E}_3,\mathrm{E}_2\}\mapsto (\mathrm{P}, a, \mathrm{E}_2), \\
&\mathsf{diag}\{\mathrm{Q},\mathrm{K},\mathrm{KJ},\mathrm{J},\mathrm{P}\}\mapsto (\mathrm{E}_2,b, \mathrm{PQ})  & \mathsf{diag}\{\mathrm{E}_2,\mathrm{E}_3, \mathrm{K},\mathrm{K},\mathrm{Q}\}\mapsto (\mathrm{Q}, c, \mathrm{E}_2)
\end{array}
\]
Then $\varphi$ is an embedding. Hence $(\hypo_3, ^{\sharp}) \in \var((A,^*)\times (B,^*))$. Therefore $\var(\hypo_3, ^{\sharp})=\var((A,^*)\times (B,^*))$. Now the result follows from Theorems~\ref{thm:balanced} and \ref{thm:B}.
\end{proof}

Let $C$  be the 25-element monoid $A_0^1\times A_0^1$.
The monoid $C$ forms an involution monoid $(C, ^*)$ with the unary
operation $^*$ given by $(x, y)^*=(y^*, x^*)$. Clearly, $(C, ^*)$ is isomorphic to the involution monoid generated by $\mathrm{E}_{6}$, $\mathsf{diag}\{\mathrm{J},\mathrm{E}_3\}$, $\mathsf{diag}\{\mathrm{K},\mathrm{E}_3\}, \mathsf{diag}\{\mathrm{E}_3, \mathrm{J}\}, \mathsf{diag}\{\mathrm{E}_3, \mathrm{K}\}$ under the skew transposition.  It is routine to verify that $(C, ^*)$ is a homomorphic image of $(\hypo_4,^\sharp)$ under the map given by $[\varepsilon]_{\hypo_4}\mapsto(1,1)$,  $[1]_{\hypo_4}\mapsto(a,1)$, $[2]_{\hypo_4}\mapsto(b,1)$, $[3]_{\hypo_4}\mapsto(1,a)$ and $[4]_{\hypo_4}\mapsto(1,b)$. Note that $(C, ^*)$ is also a homomorphic image of $(H,^\sharp)$ under the map given by $1\mapsto(1,1)$,  $a\mapsto(a,1)$, $b\mapsto(b,1)$, $c\mapsto(1,a)$ and $d\mapsto(1,b)$.

\begin{theorem}\label{thm:C}
A word identity $\bfu\approx \bfv$ holds in $(C, ^*)$ if and only if
\begin{enumerate}[\rm(i)]
  \item $\con(\bfu)=\con(\bfv), \ml(\bfu)=\ml(\bfv), \lin(\bfu)=\lin(\bfv)$;
  \item for any $x, y \in \con(\bfu)$, $x\prec_{\bfu} y$ if and only if $x\prec_{\bfv} y$.
\end{enumerate}
\end{theorem}

\begin{proof}
Let $\bfu\approx \bfv$ be a word identity satisfied by $(C, ^*)$.
Clearly, the involution submonoid of $(C, ^*)$ generated by $(a,a)$ and $(b,b)$ is isomorphic to $(A_0^1, ^*)$.  It follows from Theorem \ref{thm:A01} that $\con(\bfu)=\con(\bfv), \lin(\bfu)=\lin(\bfv)$. Suppose that $\ml(\bfu)\neq\ml(\bfv)$.  Then there exists $x\in \mix(\bfu)$ such that $\occ(x,\bfu)=1$ but $\occ(x,\bfv)\geq 2$. Let $\varphi_1$ be the homomorphism from $(\mathcal{X}\cup \mathcal{X}^*)^{+}$ to $(C,^*)$ that maps $x$ to $(1, ab)$ and any other variable to $(1,1)$. Then $\varphi_1(\bfu)\in \{(ab,ab), (ba,ab)\}, \varphi_1(\bfv)\in \{(ab,ba), (ba,ba)\}$, which implies that $\varphi_1(\bfu) \ne \varphi_1(\bfv)$, a contradiction. Therefore the condition (i) holds.

Suppose that there exist $x, y \in \con(\bfu)$ such that $x\prec_{\bfu} y$ but $x\not\prec_{\bfv} y$. Then $\bfu[x, y] \in \{x, x^*,y^*\}^{+} \cdot \{y, x^*, y^*\}^{+}$ while some $x$ occurs after the first $y$ in $\bfv$.  Let $\varphi_2$ be the homomorphism from $(\mathcal{X}\cup \mathcal{X}^*)^{+}$ to $(C,^*)$ that maps $x$ to $(a, 1)$, $y$ to $(b, 1)$ and any other variable to $(1,1)$. Then $\varphi_2(\bfu)\in\{(ab, 1), (ab, a), (ab, b), (ab, ab), (ab, ba)\}$ and $\varphi_2(\bfv)\in\{(ba, 1), (ba, a), (ba, b), (ba, ab), (ba, ba)\}$, which implies that $\varphi_2(\bfu) \ne \varphi_2(\bfv)$, a contradiction. Therefore the condition (ii) holds.

Conversely, let $\bfu\approx \bfv$ be any word identity  satisfying conditions (i) and (ii) and $\phi$ be any homomorphism from $(\mathcal{X}\cup \mathcal{X}^*)^{+}$ to $(C,^*)$.
A word $\bfu$ is a \textit{scattered subword} of $\bfw$ if there exist $\bfu_1,\dots,\bfu_n, \bfw_0, \bfw_1,\dots, \bfw_n\in F_{\mathsf{inv}}^{\varepsilon}(\mathcal{X})$ such that $\bfu=\bfu_1\cdots\bfu_n$ and $\bfw=\bfw_0\bfu_1\bfw_1\cdots\bfw_{n-1}\bfu_n\bfw_n$. The set of all scattered subwords of $\bfw$ of length $k$
is denoted by $\mathsf{sca}_k(\bfw)$. Recall that $(C, ^*)$ is isomorphic to the involution monoid generated by $\mathrm{E}_{6}$, $\mathsf{diag}\{\mathrm{J},\mathrm{E}_3\}$, $\mathsf{diag}\{\mathrm{K},\mathrm{E}_3\}, \mathsf{diag}\{\mathrm{E}_3, \mathrm{J}\}, \mathsf{diag}\{\mathrm{E}_3, \mathrm{K}\}$ under the skew transposition. It follows from \cite[Corollary~3.3]{MP2019} that to show $\phi(\bfu)=\phi(\bfv)$, it suffices to show $\mathsf{sca}_k(\bfu) =\mathsf{sca}_k(\bfv)$ for each $k=1,2,3,4,5$. Further, since each matrix  in the involution semigroup generated by $\mathrm{E}_{6}$, $\mathsf{diag}\{\mathrm{J},\mathrm{E}_3\}$, $\mathsf{diag}\{\mathrm{K},\mathrm{E}_3\}, \mathsf{diag}\{\mathrm{E}_3, \mathrm{J}\}, \mathsf{diag}\{\mathrm{E}_3, \mathrm{K}\}$ is a diagonal block matrix, we only need to show that $\mathsf{sca}_k(\bfu) =\mathsf{sca}_k(\bfv)$ for each $k=1,2$. Clearly $\mathsf{sca}_1(\bfu) =\mathsf{sca}_1(\bfv)$ by $\con(\bfu)=\con(\bfv)$.  Note that $\mathsf{sca}_2(\bfu)=\{x^2, xy|~\mbox{for any}~ x,y \in \con(\bfu)\}$. If $x^2\in \mathsf{sca}_2(\bfu)$ for some $x$, then $x^2\in \mathsf{sca}_2(\bfv)$ by (i); if $xy\in \mathsf{sca}_2(\bfu)$, then $xy\in \mathsf{sca}_2(\bfv)$ by (i) and (ii). Hence $\mathsf{sca}_2(\bfu) \subseteq \mathsf{sca}_2(\bfv)$, and  $\mathsf{sca}_2(\bfv) \subseteq \mathsf{sca}_2(\bfu)$ by symmetry. Therefore, $\mathsf{sca}_2(\bfu) = \mathsf{sca}_2(\bfv)$, as required.
\end{proof}

\begin{theorem}\label{hypo4+}
A word identity $\bfu\approx \bfv$ holds in $(\hypo_n, ^\sharp)$ for $n\geq 4$ if and only if
\begin{enumerate}[\rm(i)]
  \item $\bfu\approx \bfv$ is balanced;
  \item for any $x, y \in \con(\bfu)$, $x\prec_{\bfu} y$ if and only if $x\prec_{\bfv} y$.
\end{enumerate}
\end{theorem}

\begin{proof}
It follows from the proof of Theorem \ref{thm:same-var} that $\var(\hypo_n, ^{\sharp})=\var(H, ^{\sharp})$ for $n\geq 4$. Clearly $(A,^*), (C,^*)\in \var(H, ^{\sharp})$, thus $(A,^*)\times (C,^*)\in \var(H, ^{\sharp})$. Note that $(A, ^*)$ is isomorphic to the involution matrix monoid generated by $\mathrm{P}, \mathrm{Q}, \mathrm{E}_2$ under the skew transposition.
Let $\varphi$ be a map from $(H,^{\sharp})$ to $(A,^*)\times (A,^*)\times (C,^*)$ given by
\begin{gather*}
1\mapsto (\mathrm{E}_2,\mathrm{E}_2, (1,1)),\;\;  a \mapsto (\mathrm{P}, \mathrm{E}_2, (a,1)), \\
b\mapsto(\mathrm{E}_2,\mathrm{P}, (b,1)),\;\; c \mapsto (\mathrm{E}_2, \mathrm{Q}, (1,a)),\;\; d \mapsto (\mathrm{Q}, \mathrm{E}_2, (1,b)).
\end{gather*}
Then $\varphi$ is an embedding.  Hence $(H, ^{\sharp})\in \var((A,^*)\times (C,^*))$. Therefore $\var(H, ^{\sharp})= \var((A,^*)\times (C,^*))$. Now the result follows from Theorems~\ref{thm:balanced} and \ref{thm:C}.
\end{proof}

\begin{remark}
It follows from Theorems~\ref{thm:A01}, \ref{thm:B} and \ref{thm:C} that the involution monoids $(A_0^1,^*)$, $(B,^*)$ and $(C,^*)$ generate different varieties. However, the monoids $A_0^1$, $B$ and $C$ generate the same variety. Clearly $A_0^1 \in \var B$. Recall that $(A_0^1, ^*)$ is isomorphic to the involution matrix monoid generated by $\mathrm{J}, \mathrm{K}, \mathrm{E}_3$ under the skew transposition and $(B, ^*)$ is isomorphic to the involution matrix monoid generated by $\mathrm{E}_{9}$, $\mathsf{diag}\{\mathrm{J},\mathrm{J},\mathrm{E}_3\}, \mathsf{diag}\{\mathrm{K},\mathrm{KJ},\mathrm{J}\}, \mathsf{diag}\{\mathrm{E}_3$, $\mathrm{K},\mathrm{K}\}$ under the skew transposition. Hence $B$ is a submonoid of $A_0^1\times A_0^1\times A_0^1$, and so $B\in \var A_0^1$. Therefore $\var A_0^1=\var B$. It is obvious that $\var A_0^1=\var C$. Therefore $\var A_0^1=\var B=\var C$. Furthermore, the word identities satisfied by $A_0^1$ or $B$ or $C$ can be characterized as follows: an identity $\bfu\approx \bfv$ holds in $A_0^1$ [resp. $B, C$] if and only if
\begin{enumerate}[\rm(i)]
  \item $\con(\bfu)=\con(\bfv), \lin(\bfu)=\lin(\bfv)$;
  \item for any $x, y \in \con(\bfu)$, $x\prec_{\bfu} y$ if and only if $x\prec_{\bfv} y$.
\end{enumerate}
This fact also has been shown in \cite[Proposition~4.2(ii) and (iv)]{Sapir15}.
\end{remark}

\begin{remark}
It follows from Theorems~\ref{thm:hypo2}, \ref{hypo3} and \ref{hypo4+} that the involution monoids $(\hypo_2,^{\sharp})$, $(\hypo_3,^{\sharp})$ and $(\hypo_n,^{\sharp})$ with $n\geq 4$ generate different varieties. However, the monoids $\hypo_n$ with $n\geq 2$ generate the same variety \cite[Theorem~3.7]{CMR21a}. By Theorems~\ref{thm:hypo2}, \ref{hypo3} and \ref{hypo4+}, the word identities satisfied by $\hypo_n$ with $n\geq 2$ can be characterized as follows: an identity $\bfu\approx \bfv$ holds in $\hypo_n$ with $n\geq 2$ if and only if
\begin{enumerate}[\rm(i)]
  \item $\bfu\approx \bfv$ is balanced;
  \item for any $x, y \in \con(\bfu)$, $x\prec_{\bfu} y$ if and only if $x\prec_{\bfv} y$.
\end{enumerate}
This fact also has been shown in \cite[Theorem~4.1]{CMR21a}.
\end{remark}

\section{Finite basis problem for $(\hypo_n,^\sharp)$} \label{sec:FBP4} %
In this section, the finite basis problem for the involution monoid $(\hypo_n,^\sharp)$ for each finite $n$ is solved.
Clearly, $(\hypo_1,^\sharp)$ is finitely based since its involution is trivial and it is commutative.

To show that $(\hypo_2,^\sharp)$ and $(\hypo_3,^\sharp)$ are non-finitely based, a sufficient condition under which an involution monoid is non-finitely based is needed.
For each $k\geq 2$, define
\[
\mathcal{X}_k = \{ x_1, x_2, \ldots, x_k, \ x_1^*, x_2^*, \ldots, x_k^*, \ y_1, y_2, \ldots, y_k, \ y_1^*, y_2^*, \ldots, y_k^* \}.
\]
Then let~$\mathsf{P}_k$ denote the set of all words $\bfw \in \mathcal{X}_k^+$ such that
\[
\begin{array}{cccccll}
&        & \makebox[0.4in][r]{$\{ {_\infty}y_1, {_\infty}y_2, \ldots, {_\infty}y_k, \ {_\infty}x_1 \}$} & \precw & {{_1}y_1^*} \\[0.04in]
& \precw & {_\infty}x_2     & \precw & {_1}y_2^*\\[0.04in]
& \precw & {_\infty}x_3     & \precw & {_1}y_3^* \\[0.04in]
&        & \vdots           &        & \vdots \\[0.04in]
& \precw & {_\infty}x_{k-1} & \precw & {_1}y_{k-1}^* \\[0.04in]
& \precw & {_\infty}x_k     & \precw & \makebox[0.3in][l]{$\{ {_1}y_k^*, \ {_1}x_1^*, {_1}x_2^*, \ldots, {_1}x_k^* \}$,}
\end{array}
\]
and dually, let~$\mathsf{Q}_k$ denote the set of all words $\bfw \in \mathcal{X}_k^+$ such that
\[
\begin{array}{cccccll}
&        & \makebox[0.4in][r]{$\{ {_\infty}x_1, {_\infty}x_2, \ldots, {_\infty}x_k, \ {_\infty}y_k \}$} & \precw & {_1}x_k^* \\[0.04in]
& \precw & {_\infty}y_{k-1} & \precw & {_1}x_{k-1}^* \\[0.04in]
& \precw & {_\infty}y_{k-2} & \precw & {_1}x_{k-2}^* \\[0.04in]
&        & \vdots           &        & \vdots \\[0.04in]
& \precw & {_\infty}y_2     & \precw & {_1}x_2^* \\[0.04in]
& \precw & {_\infty}y_1     & \precw & \makebox[0.3in][l]{$\{ {_1}x_1^*, \ {_1}y_1^*, {_1}y_2^*, \ldots, {_1}y_k^* \}$.}
\end{array}
\]

\begin{lemma}[{\cite[Theorem~3.1]{GZL-A01}}]\label{lem:suffc}
Suppose that $(S, ^{*})$ is any involution monoid such that
\begin{enumerate}[{\rm(I)}]
  \item for each $k \geq 2$, there exist $\bfp_k \in \mathsf{P}_k$ and $\bfq_k \in \mathsf{Q}_k$ such that $(S, ^{*})$ satisfies
the identity $\bfp_k \approx \bfq_k$;
  \item  if $(S, ^{*})$ satisfies some word identity $\bfu \approx \bfv$ with
\[
{_\infty}x \prec_{\bfu} {_1}x^* \prec_{\bfu} {_\infty}y \prec_{\bfu} {_1}y^*,
\]
then either
\[
{_\infty}x \prec_{\bfv} {_1}x^* \prec_{\bfv} {_\infty}y \prec_{\bfv} {_1}y^*\quad \text{or} \quad  {_\infty}y \prec_{\bfv} {_1}y^* \prec_{\bfv} {_\infty}x\prec_{\bfv} {_1}x^*.
\]
\end{enumerate}
Then $(S, ^{*})$ is non-finitely based.
\end{lemma}

\begin{theorem}\label{thm: A01 NFB}
Let $\bfp_k\in\mathsf{P}_k, \bfq_k\in \mathsf{Q}_k$ for each $k\geq 2$. Then
every involution monoid $(S, ^{*})$ such that
\[
(A_0^1,^{*}) \in \var(S, ^{*})\subseteq  \var\{\mathbf{p}_k \approx \mathbf{q}_k \,|\, k\geq 2 ~\mbox{and}~ k\in \mathbb{N}\}
\]
is non-finitely based.
\end{theorem}

\begin{proof}
It follows from Theorem~\ref{thm:A01} or see \cite[Lemma~4.4]{GZL-A01} that $(A_0^1,^{*})$ satisfies the condition (II) of Lemma~\ref{lem:suffc}.
Since $(A_0^1,^{*}) \in \var(S, ^{*})$, it follows that $(S, ^{*})$ satisfies the condition (II) of Lemma~\ref{lem:suffc}.
And for each $k\geq 2$, the identity $\mathbf{p}_k \approx \mathbf{q}_k$ satisfies the condition (I) of Lemma~\ref{lem:suffc}.
Therefore $(S, ^{*})$ is  non-finitely based by Lemma~\ref{lem:suffc}.
\end{proof}

\begin{theorem}\label{thm:hypo23NFB}
Any variety in the interval $[\var(A_0^1,^{*}), \var(\hypo_3,^\sharp)]$ is non-finitely based. Consequently, $\var(B,^{*})$, $\var(\hypo_2,^\sharp)$ and $\var(\hypo_3,^\sharp)$ are non-finitely based.
\end{theorem}

\begin{proof}
For each $k\geq 2$, let
\begin{align*}
  \bfp_k=x_1x_2\cdots x_ky_1y_2\cdots y_k\cdot\bfl_k\cdot x_1^*x_2^*\cdots x_k^* y_1^*y_2^*\cdots y_k^*,\\
  \bfq_k=x_1x_2\cdots x_ky_1y_2\cdots y_k\cdot\bfr_k\cdot x_1^*x_2^*\cdots x_k^* y_1^*y_2^*\cdots y_k^*,
\end{align*}
where $\bfl_k=\lfloor \bfr^*_k \rfloor=y_1y_2\cdots y_k\cdot(x_1y_1^*\cdot x_2y_2^*\cdots x_ky_k^*)\cdot x_1^*x_2^*\cdots x_k^*$. It is routine to show that  $\bfp_k\in \mathsf{P}_k$ and $\bfq_k\in \mathsf{Q}_k$.
It follows from Theorem \ref{hypo3} that $(\hypo_3,^\sharp)$ satisfies the identity $\bfp_k\approx \bfq_k$ for each $k\geq 2$ and so each subvariety of $\var(\hypo_3,^\sharp)$ also satisfies the identity $\bfp_k\approx \bfq_k$. Therefore for any variety $\var(S, ^{*})\in [\var(A_0^1,^{*}), \var(\hypo_3,^\sharp)]$, it follows from Theorem~\ref{thm: A01 NFB} that $\var(S, ^{*})$ is non-finitely based. Note that $\var(A_0^1,^{*})$ $\subset \var(\hypo_2,^\sharp) \subset \var(\hypo_3,^\sharp)$ and $\var(A_0^1,^{*})$ $\subset \var(B,^*) \subset \var(\hypo_3,^\sharp)$. Therefore, $\var(B,^{*})$, $\var(\hypo_2,^\sharp)$ and $\var(\hypo_3,^\sharp)$ are non-finitely based.
\end{proof}

A pair $\{{_i}c, {_j}d\} \subseteq \ocs(\bfu)$ is \textit{unstable} in a balanced identity $\bfu \approx \bfv$ if ${_i}c \prec_\bfu {_j}d$ but ${_j}d \prec_\bfv {_i}c$. The set of all unstable pairs in $\bfu \approx \bfv$ will be referred as $\chaos(\bfu\approx \bfv)$. A pair $\{{_i}c, {_j}d\} \subseteq \ocs(\bfu)$ is \textit{critical} in a balanced identity $\bfu \approx \bfv$ if $\{{_i}c, {_j}d\}$ is adjacent in $\bfu$ and unstable in $\bfu \approx \bfv$.

\begin{lemma}\label{lem:unstable}
If $\{{_i}c, {_j}d\} \subseteq \ocs(\bfu)$ is unstable in a balanced word
identity $\bfu \approx \bfv$ and ${_i}c \prec_\bfu {_j}d$ then there exist a pair $\{{_s}p, {_t}q\} \subseteq \ocs(\bfu)$ such that ${_s}p \prec_\bfu {_t}q$ and $\{{_s}p, {_t}q\}$ is critical in $\bfu \approx \bfv$.
\end{lemma}
\begin{proof}
For every non-trivial balanced plain word identity, the result holds by \cite[Lemma~3.2]{Sapir2000}. Clearly, if $\bfu\approx \bfv$ is any balanced word identity, then by the proof of \cite[Lemma~3.2]{Sapir2000}, it is easy to see that the result still holds.
\end{proof}

\begin{theorem}\label{thm:C's basis}
The identities \eqref{id: inv} and
\begin{align*}
xz\,x\,tx &\approx xz\,tx, \tag{5.1}\label{e1}\\
xy\,zxty &\approx yx\,zxty, \tag{5.2}\label{e2}\\
xzyt\,xy &\approx xzyt\,yx \tag{5.3}\label{e3}
\end{align*}
constitute an identity basis for $(C,^*)$.
\end{theorem}

\begin{proof}
It follows from Theorem~\ref{thm:C} that $(C,^*)$ satisfies the identities \eqref{id: inv} and \eqref{e1}--\eqref{e3}. Note that the identities \eqref{id: inv} can be used to convert any nonempty term into some unique word. It suffices to show that each non-trivial word identity $\bfu \approx \bfv$ satisfied by $(C, \op)$ can be deduced from \eqref{e1}--\eqref{e3}. If $\bfu\approx \bfv$ is not balanced, then it follows from Theorem~\ref{thm:C}(i) that there exists some $x\in \con(\bfu)$ such that $\occ(x,\bfu), \occ(x,\bfv)\geq 2$ and $\occ(x,\bfu)\neq \occ(x, \bfv)$.  By using the identity \eqref{e1}, we can delete all of the non-first and non-last occurrences of $x$ in both $\bfu$ and $\bfv$ such that $x$ occurs exactly twice. By repeating this process, $\bfu \approx \bfv$ can be converted into a balanced word identity. Therefore we may assume that $\bfu \approx \bfv$ is a non-trivial balanced word identity.
It follows from Theorem~\ref{thm:C}(ii) that $\bfu[x, \lin(\bfu), \ml(\bfu)] = \bfv[x, \lin(\bfu), \ml(\bfu)]$ for any $x\in \con(\bfu)$ with $\occ(x, \bfu)=2$.
Since $\bfu \approx \bfv$ is non-trivial, it follows from Lemma~\ref{lem:unstable} that $\bfu \approx \bfv$ contains a critical pair $\{{_i}c, {_j}d\} \subseteq \ocs(\bfu)$. By Theorem~\ref{thm:C}, $\occ(c, \bfu)=\occ(d, \bfu)= 2$ and $(i,j)\notin \{(1, 2), (2, 1)\}$, and so $(i,j)\in \{(1, 1), (2, 2)\}$.
Therefore, one can swap ${_i}c$ and ${_j}d$ in $\bfu$ by using identities \eqref{e2} and \eqref{e3} and obtain a new word $\bfu_1$ such that $\chaos(\bfu_1 \approx \bfv_1=\bfv) \cup \{{_i}c, {_j}d\}= \chaos(\bfu \approx \bfv)$. If $\chaos(\bfu_1 \approx \bfv_1)=\emptyset$, then $\bfu_1=\bfv_1$ and so $\bfu \approx \bfv$ can be deduced from \eqref{e2} and \eqref{e3}. Otherwise, $\bfu_1 \approx \bfv_1$ is a non-trivial balanced word identity satisfying $\bfu_1[x, \lin(\bfu), \ml(\bfu)] = \bfv_1[x, \lin(\bfu), \ml(\bfu)]$ for any $x\in \con(\bfu)$ with $\occ(x, \bfu)\geq 2$. Then by using identities \eqref{e2} and \eqref{e3}, we can get the identity $\bfu_2 \approx \bfv_2$ such that $\chaos(\bfu_2 \approx \bfv_2)\subset \chaos(\bfu_1 \approx \bfv_1)$ and $|\chaos(\bfu_2 \approx \bfv_2)|+1 =|\chaos(\bfu_1 \approx \bfv_1)|$. By repeating this procedure $|\chaos(\bfu \approx \bfv)|=k$ times, we can get the identity $\bfu_k \approx \bfv_k$ such that $\chaos(\bfu_k \approx \bfv_k)=\emptyset$. Therefore, $\bfu \approx \bfv$ can be deduced from \eqref{e2} and \eqref{e3}.
\end{proof}

\begin{remark}
Note that there exists another involution operation on $C$ which can be defined by $(x,y)^{\circledast}=(x^*,y^*)$. Clearly involution semigroups $(C, ^*)$ and $(C, ^{\circledast})$ are not isomorphic and anti-isomorphic although they have the same semigroup reduct. It is easy to see that $\var(C, ^*)$ is finitely based by Theorem~\ref{thm:C's basis} while $\var(C, ^{\circledast})=\var(A_0^1,^*)$ is non-finitely based by \cite{GZL-A01}.
\end{remark}

\begin{theorem}\label{thm:hypo4's basis}
The identities \eqref{id: inv}, \eqref{e2}, \eqref{e3} and
\begin{align}
xz\,xy\,tx &\approx xz\,yx\,tx, \tag{5.4}\label{M}
\end{align}
constitute an identity basis for $(\hypo_n,^\sharp)$ with $n\geq 4$.
\end{theorem}

\begin{proof}
It follows from Theorem~\ref{hypo4+} that $(\hypo_4,^\sharp)$ satisfies the identities \eqref{id: inv}, \eqref{e2}, \eqref{e3} and \eqref{M}.
Note that the identities \eqref{id: inv} can be used to convert any nonempty term into some unique word. It suffices to show that each non-trivial word identity $\bfu \approx \bfv$ satisfied  by $(\hypo_4,^\sharp)$ can be deduced from \eqref{e2}, \eqref{e3} and \eqref{M}.
If there exists some $x\in \con(\bfu)$ such that $\bfu[x, \lin(\bfu), \ml(\bfu)] \ne  \bfv[x, \lin(\bfu), \ml(\bfu)]$, then $\occ(x,\bfu)> 2$ by Theorem~\ref{hypo4+}(ii) and there exists $t\in \lin(\bfu)\cup\ml(\bfu)$ such that ${_i}x \prec_{\bfu} {_1}t$ while ${_1}t \prec_{\bfv}{_i}x$ for some $1\leq i \leq \occ(x,\bfu)$. It follows from Theorem~\ref{hypo4+}(ii) that $i \not\in \{1, \occ(x,\bfu)\}$. Therefore by \eqref{M}, $\bfu\approx \bfv$ can be converted into an identity $\bfu'\approx \bfv'$ such that $\bfu'[x, \lin(\bfu'), \ml(\bfu')] = \bfv'[x, \lin(\bfv'),\ml(\bfu')]$ for any $x\in \con(\bfu)$.
Now the case is the same as Theorem~\ref{thm:C's basis}, and so $\bfu'\approx \bfv'$ can be derived from
\eqref{e2} and \eqref{e3}, as required.
\end{proof}

Note that the proof of \cite[Theorem 4.8]{CMR21a} can also be adapted to give
an alternative proof for Theorem \ref{thm:hypo4's basis} immediately, as it uses the same characterization of
identities.

An immediate consequence of  Theorems~\ref{thm:hypo23NFB} and \ref{thm:hypo4's basis} is the following.

\begin{corollary}
The involution monoid $(\hypo_n,^\sharp)$ has finite axiomatic rank if and only if $n\neq 2, 3$.
\end{corollary}

An involution semigroup $(S, \op)$ is \textit{twisted} if its variety $\var(S,\op)$ contains the involution semilattice $(S\ell_3, \op)$ where
\[
S\ell_3=\langle e, f~|~ e^2=e, f^2=f, ef=fe=0 \rangle,
\]
and $^*$ interchanges $e$ and $f$.
By \cite[Lemmas~7(ii), 10 and 12]{Lee17a}, we have the following result.

\begin{lemma}\label{lem:twisted}
Let $(S, \op)$ be any twisted involution semigroup. If $(S, \op)$ is finitely based, then the finite basis of $(S, \op)$ can be converted into the form
\[
\Sigma=\Sigma_{\mathsf{mix}}\cup \Sigma_{\mathsf{P}},
\]
where $\Sigma_{\mathsf{mix}}$ contains only mixed identities, and $\Sigma_{\mathsf{P}}$ contains only plain identities and is exactly a finite basis for the reduct $S$.
\end{lemma}

Let $\varphi$ be the mapping from $(\hypo_2\backslash[\varepsilon]_{\hypo_2},^{\sharp})$ to $(S\ell_3,\op)$
which maps $[1^m]_{\hypo_2}$ to $e$ and $[2^m]_{\hypo_2}$ to $f$ for all $m\geq 1$, and otherwise to $0$.
Then $\varphi$ is an epimorphism from $(\hypo_2\backslash[\varepsilon]_{\hypo_2},^{\sharp})$ to $(S\ell_3,\op)$.  Thus $(S\ell_3,\op) \in \var(\hypo_2,^{\sharp})$, and so $(\hypo_n,^{\sharp})$ is twisted for all finite $n\geq 2$. By Lemma \ref{lem:twisted} and \cite[Proposition~4.13]{CMR21a}, the identity basis for $(\hypo_n,^\sharp)$ with $n\geq 4$ contains the identities \eqref{e2}--\eqref{M} which constitute an identity basis for $\hypo_n$ when $n\geq 2$. Since the axiomatic rank of $\hypo_n$ with $n\geq 2$ is $4$ \cite[Corollary~4.12]{CMR21a}, we have the following result.

\begin{theorem}
The axiomatic rank of $(\hypo_n,^\sharp)$ with $n\geq 4$ is $4$.
\end{theorem}

By \cite[Theorem~3.2]{Jack00} and \cite[Theorems 4.1 and 4.8]{CMR21a}, it is easy to see that the variety $\var(\hypo_n)$ for each $n\geq 2$ contains continuum many subvarieties.
Now we consider the number of subvarieties of $\var(\hypo_n,^{\sharp})$.
Recall that a word $\bfu$ is an \textit{isoterm} for an involution monoid if it does not satisfy any non-trivial word
identity of the form $\bfu \approx \bfv$.

\begin{lemma}[{\cite[Theorem~3.6]{GZL-variety}}]\label{lem:subvariety}
 Let $(M, ^*)$ be any involution monoid with isoterms $xx^*yy^*$ and $xyy^*x^*$.
Then the variety $\var(M,^*)$ contains continuum many subvarieties.
\end{lemma}

\begin{theorem}\label{th:subvarieties}
Each of varieties $\var(B,^*), \var(C,^*)$ and $\var(\hypo_n,^{\sharp})$ with $n\geq 3$ contains continuum many subvarieties.
\end{theorem}

\begin{proof}
By Theorems~\ref{thm:B}--\ref{hypo4+}, it is routine to show that both the words $xx^*yy^*$ and $xyy^*x^*$ are isoterms for involution monoids $(B,^*), (C,^*)$ and $(\hypo_n,^{\sharp})$ with $n\geq 3$. Now the results follows from Lemma~\ref{lem:subvariety}.
\end{proof}

Since $xx^*yy^*$ and $xyy^*x^*$ are not isoterms for $(\hypo_2,^{\sharp})$, Lemma~\ref{lem:subvariety} can not be applied to $\var(\hypo_2,^{\sharp})$. The number of subvarieties of $\var(\hypo_2,^{\sharp})$ is still unknown.

\begin{remark}
In \cite[Question~1.5]{Lee19+}, Lee asked whether there exists a finitely based finite involution monoid such that its
variety contains continuum many subvarieties. And in the same paper, he gave an example of such an involution monoid of order $31$. By Theorems~\ref{thm:C's basis} and \ref{th:subvarieties}, the involution monoid $(C,^*)$ of order $25$ is finitely based which contains continuum many subvarieties. Hence $(C,^*)$ is another example to answer the question of \cite[Question~1.5]{Lee19+}.
\end{remark}

\section{Recognizing identities of $(\hypo_n,^\sharp)$ in polynomial time}\label{sec:CHECK-ID}%

In this section, it is shown that the identity checking problem of $(\hypo_n,^\sharp)$
belong to the complexity class $\mathsf{P}$.

\begin{theorem}\label{thm:checkABC}
The decision problems ${\textsc{Check-Id}}(A_0^1,^*), {\textsc{Check-Id}}(B,^*)$ and  ${\textsc{Check-Id}}(C,^*)$
belong to the complexity class $\mathsf{P}$.
\end{theorem}

\begin{proof}
For ${\textsc{Check-Id}}(A_0^1,^*)$, it suffices to show that, given any word identity $\bfu\approx \bfv$, one can check whether or not the words $\bfu$ and  $\bfv$ satisfy conditions of Theorem \ref{thm:A01} in polynomial time of the sum of the lengths of $\bfu$ and $\bfv$ time. For this, it suffices to exhibit algorithms that, given a word $\bfw=z_1z_2\cdots z_n$, calculate $\con(\bfw), \mix(\bfw)$ and $\lin(\bfw)$ and check whether or not $\{x, y\}\prec_\bfw\{x^*, y^*\}$ when $x, y \in \mix(\bfw)$, $x \prec_\bfw\{x^*, y\}$ or $\{x, y\}\prec_\bfw x^*$ when $x \in \mix(\bfw), y\not \in \mix(\bfw)$, and $x \prec_\bfw y$ when $x, y \not \in \mix(\bfw)$ for any $x, y \in \con(\bfw)$.

To calculate $\con(\bfw)$, we initialize $\overrightarrow{\con}(\bfw)=\emptyset$ and then scan the word $\bfw$ variable-by-variable from left to right. Each time when we read a variable of $\bfw$, we check whether the variable occurs in $\overrightarrow{\con}(\bfw)$, and if it does not occur, then we append the variable to $\overrightarrow{\con}(\bfw)$. Then we pass to the next variable if it exists or stop if the current variable is the last variable of $\bfw$. Clearly, at the end of the process, $\overrightarrow{\con}(\bfw)$ contains all variables that occur in $\bfw$ and so $\overrightarrow{\con}(\bfw)=\con(\bfw)$. The algorithm makes $n$ steps and on each step it compares with the current set $\overrightarrow{\con}(\bfw)$ whose cardinal number does not exceed $|\con(\bfw)|$. Hence, the time spent is linear in $|\con(\bfw)|n$.
For $\mix(\bfw)$, for any $x\in \con(\bfw)$, we only need to check whether the variable $x^*$ occurs in $\con(\bfw)$ or not. Clearly the time spent is linear in $|\con(\bfw)|(|\con(\bfw)|-1)$.
For $\lin(\bfw)$, for any $x\in \con(\bfw)\backslash\mix(\bfw)$, we only need to count the times of the variable $x$ occurring in $\bfw$. Clearly the time spent is linear in $n(|\con(\bfw)|-|\mix(\bfw)|)$. Therefore, checking whether the word identity $\bfu\approx \bfv$ satisfies the condition (i) of Theorem \ref{thm:A01} can be completed in polynomial time.

For any $x\in \con(\bfw)$, define $f(x)$ and $\ell(x)$ which are used to record the positions of ${_1}x$ and ${_\infty}x$ respectively. We scan the word $\bfw$ variable-by-variable from left to right. Each time we check whether $z_i=x$. If $z_{i}=x$, then we assign $i$ to $f(x)$ and stop the algorithm, otherwise we pass to the next variable. Hence, the time spent to obtain $f(x)$ is linear in $n$. Similarly, we scan the word $\bfw$ variable-by-variable from right to left to obtain $\ell(x)$. Clearly, the time spent to obtain $\ell(x)$ is linear in $n$.

Let $x, y \in \con(\bfw)$. If $x= y^*\in \mix(\bfw)$, then we locate the positions of ${_\infty}x$ and ${_1}x^*$ and check whether  $\ell(x) <f(x^*)$. There are $|\mix(\bfw)|$ such $x$. Hence, the time spent is linear in $|\mix(\bfw)|(2n+1)$.
If $x\neq y^*\in \mix(\bfw)$, then we locate the positions of ${_\infty}x, {_\infty}y$ and ${_1}x^*, {_1}y^*$ and check whether  $\ell(x)<f(x^*), \ell(x)<f(y^*)$, $\ell(y)<f(x^*), \ell(y)<f(y^*)$. There are $\binom{|\mix(\bfw)|}{2}- \frac{|\mix(\bfw)|}{2}$ pairs of such $x,y$. Hence, the time spent is linear in $(\binom{|\mix(\bfw)|}{2}- \frac{|\mix(\bfw)|}{2})(4n+4)$.
If $x\in \mix(\bfw), y\not\in \mix(\bfw)$, then we locate the positions of ${_\infty}x, {_1}x^*$ and ${_\infty}y, {_1}y$ and check whether $\ell(x)< f(x^*), \ell(x)<f(y)$ or $\ell(x)<f(x^*), \ell(y)<f(x^*)$.  There are $|\mix(\bfw)|(|\con(\bfw)|-|\mix(\bfw)|)$ pairs of such $x, y$. Hence, the time spent is linear in $|\mix(\bfw)|(|\con(\bfw)|-|\mix(\bfw)|)(4n+4)$.
If $x, y\not\in \mix(\bfw)$,  then we locate the positions of ${_\infty}x, {_1}y$ and ${_\infty}y, {_1}x$ and check whether $\ell(x)<f(y)$ or $\ell(y)<f(x)$.  There are $\binom{|\con(\bfw)|-|\mix(\bfw)|}{2}$ pairs of such $x, y$. Hence, the time spent is linear in $\binom{|\con(\bfw)|-|\mix(\bfw)|}{2}(4n+2)$. Therefore, checking whether $\bfu\approx \bfv$ satisfies the condition (ii) of Theorem \ref{thm:A01} can be completed in polynomial time.

Therefore the decision problem ${\textsc{Check-Id}}(A_0^1,^*)$ belongs to the complexity class $\mathsf{P}$. By a similar argument, the decision problems ${\textsc{Check-Id}}(B,^*)$ and  ${\textsc{Check-Id}}(C,^*)$ also belong to the complexity class $\mathsf{P}$.
\end{proof}

\begin{theorem}\label{thm:checkhypo}
The decision problem ${\textsc{Check-Id}}(\hypo_n,^\sharp)$ for each finite $n$
belong to the complexity class $\mathsf{P}$.
\end{theorem}

\begin{proof}
By Theorems \ref{thm:hypo2}, \ref{hypo3}, \ref{hypo4+} and \ref{thm:checkABC}, it suffices to show that, given any word identity $\bfu\approx \bfv$, there is an algorithm to check the identity $\bfu \approx \bfv$ is balanced  in polynomial time. Then we only need to check whether  $\con(\bfu) =\con(\bfv)$ and $\occ(x, \bfu)=  \occ(x, \bfv)$ for any $x\in \con(\bfu) =\con(\bfv)$. Clearly, this can be done in polynomial time.
\end{proof}

Cain et al. give a complete characterization of the identities satisfied by $\hypo_n$ for each $n\geq 2$ in \cite[Theorem~4.1]{CMR21a}. By a similar argument with Theorem \ref{thm:checkhypo}, it can be shown that the decision problem \textsc{Check-Id}($\hypo_n$) belongs to the complexity class $\mathsf{P}$.

\section*{Acknowledgements}
The authors are very grateful to the anonymous referee for his/her careful reading and suggestions.

\end{sloppypar}
\end{document}